\documentclass[10pt,a4paper]{article}

\usepackage{amsmath}
\usepackage{graphicx,latexsym,euscript,makeidx,color}
\usepackage{amsmath,amsfonts,amssymb,amsthm,mathrsfs}
\usepackage{enumerate}

\usepackage{amssymb}
\usepackage{amsmath}
\usepackage{amsfonts}
\usepackage{graphicx}
\usepackage{graphicx}          
\usepackage{color}

\usepackage[colorlinks,linkcolor=blue,anchorcolor=green,citecolor=red]{hyperref}

\usepackage{geometry}
\def\sqr#1#2{{\vcenter{\vbox{\hrule height.#2pt
              \hbox{\vrule width.#2pt height#1pt \kern#1pt \vrule width.#2pt}
              \hrule height.#2pt}}}}
\def\signed #1{{\unskip\nobreak\hfil\penalty50
              \hskip2em\hbox{}\nobreak\hfil#1
              \parfillskip=0pt \finalhyphendemerits=0 \par}}
\def\endpf{\signed {$\sqr69$}}

\def\5n{\negthinspace \negthinspace \negthinspace \negthinspace \negthinspace }
\def\4n{\negthinspace \negthinspace \negthinspace \negthinspace }
\def\3n{\negthinspace \negthinspace \negthinspace }
\def\2n{\negthinspace \negthinspace }
\def\1n{\negthinspace }

\def\dbE{\mathbb{E}}     
   \def\cF{{\cal F}}

\def\ms{\medskip}                
               
\def\ds{\displaystyle}

            \def\({\Big (}
                  \def\){\Big )}
          \def\[{\Big[}
           \def\]{\Big]}

\def\bde{\begin{definition}\label}    \def\ede{\end{definition}}
\def\bt{\begin{theorem}\label}        \def\et{\end{theorem}}
\def\bc{\begin{corollary}\label}      \def\ec{\end{corollary}}
\def\bl{\begin{lemma}\label}          \def\el{\end{lemma}}
\def\bp{\begin{proposition}\label}    \def\ep{\end{proposition}}
\def\bas{\begin{assumption}\label}    \def\eas{\end{assumption}}
\def\br{\begin{remark}\label}         \def\er{\end{remark}}
\def\bex{\begin{example}\label}       \def\ex{\end{example}}
\def\ba{\begin{array}}                \def\ea{\end{array}}
\def\be{\begin{equation}}
\def\bel{\begin{equation}\label}      \def\ee{\end{equation}}
\def\bea{\begin{eqnarray*}}           \def\eea{\end{eqnarray*}}

\def\square#1{\vbox{\hrule\hbox{\vrule height#1%
     \kern#1\vrule}\hrule}}
\def\rectangle#1#2{\vbox{\hrule\hbox{\vrule height#1%
     \kern#2\vrule}\hrule}}

\font\tenbb=msbm10 \font\sevenbb=msbm7 \font\fivebb=msbm5
\newfam\bbfam
\scriptscriptfont\bbfam=\fivebb \textfont\bbfam=\tenbb
\scriptfont\bbfam=\sevenbb

\newtheorem{theorem}{\indent Theorem}[section]
\newtheorem{definition}[theorem]{\indent Definition}
\newtheorem{proposition}[theorem]{\indent Proposition}
\newtheorem{corollary}[theorem]{\indent Corollary}
\newtheorem{lemma}[theorem]{\indent Lemma}
\newtheorem{remark}[theorem]{\indent Remark}
\newtheorem{example}[theorem]{\indent Example}

\newtheorem{assumption}[theorem]{\indent Assumption}

\makeatletter
   
   \@addtoreset{equation}{section}
\makeatother

\allowdisplaybreaks[4]

\begin{document}

\title{\bf Time-Inconsistent Mean-Field Stochastic LQ Problem: Open-Loop Time-Consistent Control
\thanks{This work is supported in part by the National Natural Science Foundation of China (11471242, 612279002) and the National Key Basic
Research Program of China (973 Program) under grant 2014CB845301.}
}
\author{Yuan-Hua Ni\thanks{Department of Mathematics,
School of Science, Tianjin Polytechnic University, Tianjin 300387,
P.R. China. 
E-mail: {\tt yhni@amss.ac.cn}. 
}~~~~~Ji-Feng Zhang\thanks{Key Laboratory of Systems and
Control, Institute of Systems Science, Academy of Mathematics and
Systems Science, Chinese Academy of Sciences, Beijing 100190, P. R.
China. E-mail: {\tt jif@iss.ac.cn}.  
}~~~~~Miroslav Krstic\thanks{Department of Mechanical and Aerospace Engineering, University of California, San Diego, USA. E-mail: {\tt krstic@ucsd.edu}.}}
\maketitle

{\bf Abstract:} This paper is concerned with the open-loop time-consistent solution of time-inconsistent
mean-field stochastic linear-quadratic optimal control. Different from standard stochastic linear-quadratic problems, both the system matrices
and the weighting matrices are dependent on the initial times, and the conditional
expectations of the control and state enter quadratically into the cost functional. Such features will ruin Bellman's principle of optimality and result in the time-inconsistency of the optimal control. Based on the dynamical nature of the systems involved, a kind of open-loop time-consistent equilibrium control is investigated in this paper. It is shown that
the existence of open-loop time-consistent equilibrium control for a fixed initial
pair is equivalent to the solvability of a set of forward-backward stochastic
difference equations with stationary conditions and convexity conditions.
By decoupling the forward-backward stochastic
difference equations, necessary and sufficient conditions
in terms of linear difference equations and generalized difference
Riccati equations are given for the existence of open-loop time-consistent
equilibrium control with a fixed initial pair.
Moreover, the existence of open-loop time-consistent equilibrium control for
all the initial pairs is shown to be equivalent to the solvability of a set of coupled constrained
generalized difference Riccati equations and two sets of constrained linear
difference equations.

\ms

{\bf Key words:} Time-inconsistency, time-consistent solution, mean-field stochastic linear-quadratic optimal control, indefinite stochastic linear-quadratic optimal control

\ms


\section{Introduction}

Though not mentioned frequently, time-consistency is indeed an
essential notion in optimal control theory, which relates to
Bellman's principle of optimality. To see this, let us begin with a
standard discrete-time stochastic optimal control problem, whose
system dynamics and cost functional are given, respectively, by
\begin{eqnarray}\label{system-general}
\left\{\begin{array}{l}X_{k+1}=f_k(X_k,u_k,w_k),\\[1mm]
X_t=x\in \mathbb{R}^n, ~~k\in \mathbb{T}_t,~~t\in \mathbb{T},
\end{array}
\right.\end{eqnarray}
and
\begin{eqnarray}\label{performace-general}
J(t,x;u)=\sum_{k=t}^{N-1}\mathbb{E}\big{[}L_k(X_k,u_k)\big{]}+\mathbb{E}\big{[}h(X_{N})\big{]}.
\end{eqnarray}
Here, $\mathbb{T}_t=\{t,\cdots,N-1\}, \mathbb{T}=\{0,1,\cdots,N-1\}$, and $N$ is a positive integer; $\{X_k, k\in
\widetilde{{\mathbb{T}}}_t\}$ and $\{u_k, k\in {\mathbb{T}}_t\}$
with $\widetilde{\mathbb{T}}_t=\{t,t+1,\cdots,N\}$ are the state
process and the control process, respectively; $\{w_k,k\in
\mathbb{T}_t\}$ is a stochastic disturbance; $\mathbb{E}$ is the
operator of mathematical expectation. Without loss of generality, $f_k$, $L_k$, $k\in \mathbb{T}_t$, and $h$ are assumed to be bounded.
Let $\mathcal{U}[t,N-1]$ be a set of admissible controls. We then have the following optimal control problem.

\textbf{Problem (C).} Concerned with (\ref{system-general}),
(\ref{performace-general}) and the initial pair $(t, x)\in
\mathbb{T}\times \mathbb{R}^n$, find a $\bar{u}\in
\mathcal{U}[t,N-1]$ such that
\begin{eqnarray*}
J(t,x;\bar{u})=\inf_{u\in \mathcal{U}[t,N-1]}J(t,x;u).
\end{eqnarray*}
Any $\bar{u}\in\mathcal{U}[t,N-1]$ satisfying the above is called an
optimal control for the initial pair $(t, x)$;
$\bar{X}=\{\bar{X}_k=\bar{X}(k;t,x,\bar{u}), k\in
\widetilde{{\mathbb{T}}}_t\}$ is called the optimal trajectory corresponding to $\bar{u}$, and
$(\bar{X},\bar{u})$ is referred to as an optimal pair for the
initial pair $(t,x)$.

By Bellman's principle of optimality, if $\bar{u}$ is an
optimal control of Problem (C) for the initial pair $(t,x)$, then
for any $\tau\in \mathbb{T}_{t+1}=\{t+1,...,N-1\}$, $\bar{u}|_{\mathbb{T}_{\tau}}$
(the restriction of $\bar{u}$ on $\mathbb{T}_{\tau}=\{\tau,...,N-1\}$) is an optimal
control of Problem (C) for the initial pair $(\tau,
\bar{X}(\tau;t,x,\bar{u}))$.
This property is essential to handle optimal control problems like
Problem (C) and its continuous-time counterpart, which provides the
theoretical foundation of dynamic programming approach.
Such a phenomenon is referred to as the time-consistency of the
optimal control, which ensures that one needs only to solve an
optimal control problem for a given initial pair, and the obtained
optimal control is also optimal along the whole optimal trajectory.

However, in reality, the time-consistency fails quite often. For
instance, when the initial time or initial state enters into the
system dynamics or cost functional explicitly, or even more, the
conditional expectations of the state or control enters nonlinearly
into the cost functional, the corresponding problems are
time-inconsistent. See examples in \cite{Krusell} and \cite{Bjork}
about the hyperbolic discounting and quasi-geometric discounting.
The problem with nonlinear terms of conditional expectation in the
cost functional is called as mean-field stochastic optimal control.
In this case, the smoothing property of conditional expectation
will not be sufficient to ensure the time-consistency of the optimal
control. A well-known example of this case is the mean-variance utility \cite{Basak} \cite{Bjork}.

To handle the time-inconsistency, we have two different ways. The
first one is static formulation or pre-commitment formulation. If
one is able to commit to his/her initial policy and does not
revisit the problem in the future, then this policy can be
implemented as planned.
This approach neglects the time-inconsistency and the optimal
control is optimal only when viewed at the initial time. Though the
static formulation is of some practical and theoretical values, it has not
really addressed the time-inconsistency nor provided solution in a
dynamic sense.
Relative to this, another approach addresses
the time-inconsistency in a dynamic manner. Instead of seeking an
``optimal control", some kinds of equilibrium solutions are concerned
with. This is mainly motivated by practical applications such as in
mathematical finance and economics, and has recently attracted
considerable interest and efforts.
%
The mathematical formulation of the time-inconsistency was first
reported by \cite{Strotz}, and its qualitative analysis might be traced back to \cite{Smith}.
Following \cite{Strotz}, the works
\cite{Goldman}, \cite{Krusell}, \cite{Laibson} and \cite{Palacios} are for systems described by
difference equations or ordinary differential equations (ODEs).
%
%
Recently,
\cite{Ekeland} and \cite{Ekland-2} studied the non-exponential
discounting problems both for simple ODEs and stochastic
differential equations (SDEs), and introduced the notion of time-consistent
control. \cite{Bjork} discussed the problems of general Markovian
time-inconsistent stochastic optimal control.  \cite{Yong-1} and
\cite{Yong-0} addressed the deterministic continuous-time
linear-quadratic (LQ) optimal control using a cooperative
game approach. Different from \cite{Yong-1} and \cite{Yong-0},
\cite{Hu-jin-Zhou} studied another kind of time-consistent
equilibrium solution of a continuous-time time-inconsistent stochastic LQ problem. In \cite{Yong-2013}, the author investigated both the open-loop and the
closed-loop time-consistent solutions for the general mean-field
stochastic LQ problems, and showed that the
existence of open-loop equilibrium control and closed-loop
equilibrium strategy is ensured via the solvability of certain
sets of Riccati-type equations.
It is worth noting that all these existing results about LQ problems are focusing on the definite case. Here, by the definiteness, we mean that in the cost
functionals the state weight matrices are nonnegative definite and
the control weight matrices are positive definite. Furthermore, no necessary and sufficient condition is reported on the existence of time-consistent solutions for time-inconsistent LQ problems.

In this paper, we shall investigate a time-inconsistent mean-field
stochastic LQ optimal control problem, whose system dynamics and cost
functional are also dependent on the initial time. No definiteness constraint is required for either the state
or the control weighting matrices, and a class of open-loop
equilibrium control is studied for the considered LQ problem.
The main idea and results of this paper are as follows.
\begin{itemize}
\item After giving the definition of open-loop equilibrium pair,
we show in Theorem \ref{Theorem-Equivalentce-open-loop} that the existence of open-loop equilibrium pair for a fixed
initial pair is
equivalent to the solvability of a set of forward-backward
stochastic difference equations (FBS$\Delta$Es) with stationary conditions
and convexity conditions. Different from \cite{Li-Ni-Zhang}, the
equivalent conditions are proved based on a formula of cost
functional difference (Lemma \ref{Lemma-difference}).

\vskip 0.3cm
\item If for a fixed initial pair, Problem (LQ) admits an open-loop
equilibrium pair, then a set of constrained linear difference
equations (LDEs) is solvable, and the open-loop equilibrium control
admits a closed-loop representation (Theorem
\ref{Theorem-Necessary}). Here, the closed-loop representation is a
linear feedback of current value of the equilibrium state, whose
gains are computed via the solutions of a set of constrained LDEs (\ref{P}), a set of LDEs
(\ref{pi}) and a set of generalized difference Riccati equations
(GDREs) (\ref{T}).

\vskip 0.3cm
\item  If for a fixed initial pair $(t,x)$ Problem (LQ) admits
an open-loop equilibrium pair, then for any $(k, \zeta)$ with $k\in
\mathbb{T}_t$ and $\zeta\in L^2_{\mathcal{F}}(k; \mathbb{R}^n)$,
Problem (LQ) is point-wisely convex at $(k, \zeta)$. In this case,
we equivalently have the solvability of the constraint LDEs (\ref{P}).

\vskip 0.3cm \noindent Conversely, if a version of (\ref{P}) (with
$\mathbb{T}_t$ replaced by $\mathbb{T}$, i.e., (\ref{p-3})) is
solvable, then we can take a perturbation of the cost functional by
adding $\varepsilon \mathbb{E}\big{[}u_k^Tu_k\big{]}$ to $J(k,
\zeta; (u_k,...,u_{N-1}))$, $k\in \mathbb{T}$; and the
obtained problem is denoted as Problem (LQ)$_\varepsilon$, which admits an
open-loop equilibrium pair $(X^{\varepsilon, t,x,*}, u^{\varepsilon,
t,x,*})$ for
any initial pair $(t,x)$. Furthermore, if $\{u^{\varepsilon, t,x,*},
\varepsilon>0\}$ is bounded, then Problem (LQ) for the initial pair
$(t,x)$ will admit an open-loop equilibrium pair.

\vskip 0.3cm
\item
For any initial pair, Problem (LQ) admitting an open-loop
equilibrium pair is shown to be equivalent to that two sets of constrained LDEs
(\ref{p-3}) (\ref{pi-3}) and a set of constrained GDREs (\ref{T-3}) are
solvable. 
It is worth pointing out that if
solvable, the set of GDREs (\ref{T-3}) does not have symmetric
structure, i.e., its solution is not symmetric.

\end{itemize}

In \cite{Li-Ni-Zhang}, a simplified version of Problem (LQ) is
considered, where there are no mean-field terms in the system dynamics
and cost functional. Hence, this paper is a
continuation of \cite{Li-Ni-Zhang}. Concerned with the necessary and
sufficient condition on the existence of open-loop equilibrium
pair, \cite{Li-Ni-Zhang} just gives a counterpart of Corollary
\ref{corollary-Necessary-sufficient-2} of this paper with (\ref{system-1-in})
replaced by (\ref{system-1-in-2}). This is because in
\cite{Li-Ni-Zhang} we do not have a result similar to Lemma
\ref{Lemma-Z}, which gives the representation of the backward state
via the forward state. If the system dynamics and cost
functional are both independent of the initial time, the corresponding LQ
problem will be a dynamic version of that considered in
\cite{Ni-Zhang-Li}, where the conditional expectation operators are replaced by the expectation
operators. For details on mean-field stochastic optimal control and
related mean-field games, we refer to
\cite{Bensoussan} \cite{Elliott-Li-Ni}
\cite{Huang2007} \cite{Lasry} \cite{Ni-Zhang-Li} \cite{Yong-MFF-2011} and the references
therein.

Though the equilibrium {control} (\ref{open-loop-equilibrium-k}) is of feedback form, it is indeed an open-loop control. To clarify, the closed-loop expression  (\ref{open-loop-equilibrium-k}) is not a closed-loop/feedback equilibrium solution of Problem (LQ) at all. Instead, a closed-loop or feedback equilibrium solution of Problem (LQ) is concerned with the time-consistency of the \emph{strategy}. Here, by a strategy we mean a decision rule that a player/controller adopts to select her actions, based on available information. Therefore, mathematically, a strategy is a
measurable function of the information set. When the information set is available and substituted into the strategy, a control is obtained, which is then viewed as the {open-loop realization} of that strategy.
Due to their intrinsical difference between an open-loop control and a strategy, the open-loop equilibrium control of this paper differs clearly from the closed-loop/feedback equilibrium strategy, which is studied in \cite{Ni-Zhang-Krstic} by the authors.

The rest of this paper is organized as follows. Section
\ref{section-2} introduces the notion of open-loop equilibrium control of Problem (LQ), and presents necessary and sufficient conditions on its existence for both the case with a fixed initial pair and the case with all the initial pairs. Section
\ref{section-3} studies two special cases of Problem (LQ). Section
\ref{section-4} gives an example, and some concluding remarks are given in Section \ref{section-5}.

\section{Open-loop Time-Consistent Solution}\label{section-2}

Consider the following controlled stochastic difference equation
(S$\Delta$E)
\begin{eqnarray}\label{system-1}
\left\{\begin{array}{l}
X^t_{k+1}=\big{(}A_{t,k}X^t_k+\bar{A}_{t,k}\mathbb{E}_tX^t_k+B_{t,k}u_k+\bar{B}_{t,k}
\mathbb{E}_tu_k+f_{t,k}\big{)}\\[1mm]
\hphantom{X^t_{k+1}=}+\big{(}C_{t,k}X^t_k+\bar{C}_{t,k}
\mathbb{E}_tX^t_k+D_{t,k}u_k+\bar{D}_{t,k}\mathbb{E}_tu_k
+d_{t,k}\big{)}w_k, \\[1mm]
X^t_t=x,~~k\in  \mathbb{T}_t,~~t\in \mathbb{T},
\end{array}
\right.
\end{eqnarray}
where
$A_{t,k},\bar{A}_{t,k},C_{t,k},\bar{C}_{t,k}\in \mathbb{R}^{n\times
n}$, $B_{t,k},\bar{B}_{t,k},D_{t,k},\bar{D}_{t,k}\in
\mathbb{R}^{n\times m}$ are deterministic matrices, and $f_{t,k},
d_{t,k}\in \mathbb{R}^n$ are deterministic vectors; $\{X^t_{k}, k\in
\widetilde{{\mathbb{T}}}_t\}\triangleq X^t$ and $\{u_k, k\in
\mathbb{T}_t\}\triangleq u$ are the state process and the control
process, respectively. The noise $\{w_k, k\in
\mathbb{T}\}$ is assumed to be  a martingale difference sequence
defined on a probability space $(\Omega, \mathcal{F}, P)$ with
\begin{eqnarray}\label{w-moment}
\mathbb{E}_{k+1}[w_{k+1}]=0,~~\mathbb{E}_{k+1}[(w_{k+1})^2]=1,~k\geq 0.
\end{eqnarray}
$\mathbb{E}_t$ in (\ref{system-1}) is the conditional mathematical expectation
$\mathbb{E}[\,\cdot\,|\mathcal{F}_t]$ with respect to $\cF_t=\sigma\{w_l, l=0, 1,\cdots,t-1\}$, and $\mathcal{F}_{0}$ is understood as
$\{\emptyset, \Omega\}$.
The cost functional associated with the system (\ref{system-1}) is
\begin{eqnarray}\label{cost-1}
&&\hspace{-2em}J(t,x;u)
 =\sum_{k=t}^{N-1}\dbE_t\Big{\{}(X_k^t)^TQ_{t,k}X^t_k +(\mathbb{E}_tX_k^t)^T
 \bar{Q}_{t,k}\mathbb{E}_tX^t_k+ u_k^TR_{t,k}u_k+(\mathbb{E}_tu_k)^T \bar{R}_{t,k}
 \mathbb{E}_tu_k\nonumber\\[1mm]
&&\hspace{-2em}\hphantom{J(t,x;u)=}  +2q_{t,k}^T X^t_k
+ 2\rho_{t,k}^T u_k\Big{\}}+\mathbb{E}_t\big{[}(X_N^t)^TG_tX^t_N\big{]}
+(\mathbb{E}_tX^t_N)^T \bar{G}_t\mathbb{E}_tX^t_{N}+2\mathbb{E}_tg_t^T X^t_N,
\end{eqnarray}
where $Q_{t,k}, \bar{Q}_{t,k},R_{t,k},\bar{R}_{t,k}, k\in
\mathbb{T}_t$, $G_t, \bar{G}_t$ are deterministic symmetric matrices
of appropriate dimensions, and $q_{t,k}, \rho_{t,k},  k\in
\mathbb{T}_t, g_t$ are deterministic vectors. In (\ref{system-1}),
$x$ is in $L^2_\mathcal{F}(t; \mathbb{R}^n)$, which is a set of
random variables such that any $\xi\in L^2_\mathcal{F}(t; \mathbb{R}^n)$ is
$\mathcal{F}_{t}$-measurable and $\mathbb{E}|\xi|^2<\infty$. Let further
$L^2_\mathcal{F}(\mathbb{T}_t; \mathcal{H})$ be a set of
$\mathcal{H}$-valued processes such that for any its element $\nu=\{\nu_k, k\in \mathbb{T}_t\}$, $\nu_k$ is $\mathcal{F}_{k}$-measurable and
$\sum_{k=t}^{N-1}\mathbb{E}|\nu_k|^2<\infty$. Then, we pose the
following optimal control problem.

{\textbf{Problem (LQ)}.} \emph{ Concerned with
(\ref{system-1}), (\ref{cost-1}) and the initial pair $(t,x)$, find a
${u}^*\in L^2_\mathcal{F}(\mathbb{T}_t; \mathbb{R}^m)$, such that
\begin{eqnarray}\label{Problem-LQ}
J(t,x;{u}^*) = \inf_{u\in L^2_\mathcal{F}(\mathbb{T}_t; \mathbb{R}^m)}J(t,x;u).
\end{eqnarray}}

Due to the feature of time-inconsistency of Problem (LQ),
the notion ``optimality" should be defined in an appropriate way.
Therefore, instead of solving Problem (LQ) for a static
pre-committed optimal control, we adopt the concept of dynamic
equilibrium control, which is optimal in an ¡°infinitesimal¡± sense
and is  consistent with the dynamical nature of Problem (LQ).

\begin{definition}\label{Definition-open-loop}
Given $t\in \mathbb{T}$ and $x\in L^2_\mathcal{F}(t; \mathbb{R}^n)$, a
state-control pair $(X^{t,x,*}, u^{t,x,*})$ with $u^{t,x,*}\in
L^2_{\mathcal{F}}({\mathbb{T}}_t;\mathbb{R}^m)$ is called an
open-loop equilibrium pair of Problem (LQ) for the initial pair $(t,
x)$ if $X^{t,x,*}_t=x$, and
\begin{eqnarray}\label{open-loop-equilibrium}
J(k,X^{t,x,*}_k; u^{t,x,*}|_{\mathbb{T}_k})\leq J(k,X^{t,x,*}_k;
(u_k, u^{t,x,*}|_{\mathbb{T}_{k+1}}))
\end{eqnarray}
holds for any $k\in \mathbb{T}_t$ and any $u_k\in
L^2_\mathcal{F}(k; \mathbb{R}^m)$.
Here, $u^{t,x,*}|_{\mathbb{T}_k}$ and
$u^{t,x,*}|_{\mathbb{T}_{k+1}}$  (with $\mathbb{T}_k=\{k,...,N-1\}, \mathbb{T}_{k+1}=\{k+1,...,N-1\}$)  are the restrictions of $u^{t,x,*}$
on $\mathbb{T}_k$ and $\mathbb{T}_{k+1}$, respectively. Furthermore,
such a $u^{t,x,*}$ is called an open-loop equilibrium control for
the initial pair $(t,x)$, and $X^{t,x,*}$ is the corresponding equilibrium state.
\end{definition}

For any $u\in L^2_\mathcal{F}(\mathbb{T}_t; \mathbb{R}^m)$, the
requirement that $u_k$ is $\mathcal{F}_{k}$-measurable is parallel
to the standard statement on the admissible controls of
continuous-time stochastic optimal control; see \cite{Fleming},
\cite{Yong-Zhou} for details. Furthermore, $u \in
L^2_\mathcal{F}(\mathbb{T}_t; \mathbb{R}^m)$ can be viewed as an
open-loop control \cite{Astrom}.
Noting that
$u^{t,x,*}|_{\mathbb{T}_k}=(u^{t,x,*}_k,
u^{t,x,*}|_{\mathbb{T}_{k+1}})$, the control $(u_k,
u^{t,x,*}|_{\mathbb{T}_{k+1}})$ on the right-hand side of (\ref{open-loop-equilibrium}) differs from
$u^{t,x,*}|_{\mathbb{T}_k}$ only at time instant $k$. Intuitively, the cost functional will increase if one deviates from $u^{t,x,*}$.
Similar to \cite{Li-Ni-Zhang},
$\{u^{t,x,*}_t,...,u^{t,x,*}_{N-1}\}$ can be viewed as an equilibrium of a
multi-person game with hierarchical structure.
Hence, we call $u^{t,x,*}$  an open-loop equilibrium control.
By its definition, $u^{t,x,*}$ is time-consistent in the sense that
for any $k\in \mathbb{T}_t$, $u^{t,x,*}|_{\mathbb{T}_k}$ is an
open-loop equilibrium control for the initial pair
$(k,X^{t,x,*}_k)$. 

The following result is concerned with the difference of cost
functionals,  which is characterized via the solutions of an S$\Delta$E and
a backward stochastic difference equation (BS$\Delta$E).

\begin{lemma}\label{Lemma-difference}
Let $\zeta\in L^2_{\mathcal{F}}(k;\mathbb{R}^n)$, $u=\{u_\ell,k \in
\mathbb{T}_k\} \in L^2_{\mathcal{F}}(\mathbb{T}_k;\mathbb{R}^m)$,
$\bar{u}_k\in L^2_{\mathcal{F}} (k;\mathbb{R}^m)$ and  $\lambda \in
\mathbb{R}$. Then, we have
\begin{eqnarray}\label{appendix-A-J-0}
&&\hspace{-3em}{J}(k,\zeta; (u_k+\lambda \bar{u}_k,u|_{\mathbb{T}_{k+1}}))-{J}(k,\zeta; u)=2\lambda \Big{[}(R_{k,k}+\bar{R}_{k,k})u_k+(B_{k,k}+\bar{B}_{k,k})^T
\mathbb{E}_kZ^{k,u_k}_{k+1}+\rho_{k,k}\nonumber\\[1mm]
&&\hspace{-3em}\hphantom{{J}(k,\zeta; (u_k+\lambda \bar{u}_k,u|_{\mathbb{T}_{k+1}}))-{J}(k,\zeta; u)=}+(D_{k,k}+\bar{D}_{k,k})^T\mathbb{E}_k(Z_{k+1}^{k,u_k}w_k)
\Big{]}^T\bar{u}_k+ \lambda^2\widehat{{J}}(k,0;\bar{u}_k)
\end{eqnarray}
with (noting $Y^{k,\bar{u}_k}_k=0$)
\begin{eqnarray}\label{hat-J}
&&\hspace{-2em}\widehat{{J}}(k,0;\bar{u}_k)=\sum_{\ell=k}^{N-1}\mathbb{E}_k\Big{[}
(Y^{k,\bar{u}_k}_\ell)^TQ_{k,\ell} Y_\ell^{k,\bar{u}_k}+(\mathbb{E}_k
Y^{k,\bar{u}_k}_\ell)^T\bar{Q}_{k,\ell} \mathbb{E}_kY_\ell^{k,\bar{u}_k}\Big{]}\nonumber \\
&&\hspace{-2em}\hphantom{\widehat{{J}}(k,0;\bar{u}_k)=}+\mathbb{E}_k\big{[} {u}_k^T(R_{k,k}
+\bar{R}_{k,k}){u}_k \big{]}+\mathbb{E}_k\big{[}
(Y_N^{k,\bar{u}_k})^T G_{k} {Y}_N^{k,\bar{u}_k}\big{]}+(\mathbb{E}_kY_N^{k,\bar{u}_k})^T \bar{G}_{k} \mathbb{E}_k{Y}_N^{k,\bar{u}_k}.
\end{eqnarray}
Here, $u|_{\mathbb{T}_{k+1}}$
is the restriction of $u$ on $\mathbb{T}_{k+1}$, and $Z^{k,u_k}$,
$Y^{k,\bar{u}_k}$ are given, respectively, by the BS$\Delta$E
\begin{eqnarray}\label{X-Z-1}
\left\{\begin{array}{l}
Z_\ell^{k,u_k}=A_{k,\ell}^T\mathbb{E}_\ell Z_{\ell+1}^{k,u_k}+ \bar{A}_{k,\ell}^T
\mathbb{E}_kZ_{\ell+1}^{k,u_k}+C_{k,\ell}^T\mathbb{E}_\ell(Z_{\ell+1}^{k,u_k}w_\ell)
\\[1mm]
\hphantom{Z_\ell^{k,u_k}=}+\bar{C}_{k,\ell}^T\mathbb{E}_k(Z_{\ell+1}^{k,u_k}w_\ell)+Q_{k,\ell}X_\ell^{k,u_k}+\bar{Q}_{k,\ell}
\mathbb{E}_kX_\ell^{k,u_k}+q_{k,\ell},\\[1mm]
Z_N^{k,u_k}=G_kX^{k,u_k}_N+\bar{G}_k\mathbb{E}_kX^{k,u_k}_N+g_k,~~\ell\in
\mathbb{T}_k,
\end{array}\right.
\end{eqnarray}
and the S$\Delta$E
\begin{eqnarray}\label{system-y-k}
\left\{
\begin{array}{l}
{Y}^{k,\bar{u}_k}_{\ell+1}=A_{k,\ell}{Y}^{k,\bar{u}_k}_\ell+\bar{A}_{k,\ell}
\mathbb{E}_k{Y}^{k,\bar{u}_k}_\ell+\big{(}C_{k,\ell}{Y}^{k,\bar{u}_k}_\ell
+\bar{C}_{k,\ell}\mathbb{E}_k{Y}^{k,\bar{u}_k}_\ell\big{)}w_\ell,\\[1mm]
Y^{k,\bar{u}_k}_{\ell+1}=(B_{k,k}+\bar{B}_{k,k})\bar{u}_k+{(}D_{k,k}
+\bar{D}_{k,k}{)}\bar{u}_k w_k,\\[1mm]
{Y}^{k,\bar{u}_k}_k=0,~~~\ell\in
\mathbb{T}_{k+1},
\end{array}
\right.
\end{eqnarray}
where
\begin{eqnarray*}
\left\{
\begin{array}{l}
X^{k,u_k}_{\ell+1} = \big{(}A_{k,\ell}X^{k,u_k}_\ell+\bar{A}_{k,\ell}
\mathbb{E}_kX^{k,u_k}_\ell+B_{k,\ell}u_\ell+\bar{B}_{k,\ell}
\mathbb{E}_ku_\ell+f_{k,\ell}\big{)}\\[1mm]
\hphantom{X^t_{k+1}=}+\big{(}C_{k,\ell}X^{k,u_k}_\ell+\bar{C}_{k,\ell}
\mathbb{E}_kX^{k,u_k}_\ell+D_{k,\ell}u_\ell+\bar{D}_{k,\ell}\mathbb{E}_k u_\ell
+d_{k,\ell}\big{)}w_\ell,\\[1mm]
X^{k,u_k}_k=\zeta,~~ \ell\in\mathbb{T}_{k}.
\end{array}
\right.
\end{eqnarray*}

\end{lemma}

\emph{Proof.}  See Appendix A. \endpf

From Lemma \ref{Lemma-difference}, we have the following result,
which gives the necessary and sufficient condition to the existence of open-loop equilibrium pair for a given initial pair.

\begin{theorem}\label{Theorem-Equivalentce-open-loop}
Given $t\in \mathbb{T}$ and $x\in L^2_\mathcal{F}(t; \mathbb{R}^n)$, the
following statements are equivalent.

(i) There exists an open-loop equilibrium pair of Problem (LQ) for the initial pair $(t,x)$.

(ii) There exists a $u^{t,x,*}\in
L^2_{\mathcal{F}}(\mathbb{T}_t;\mathbb{R}^m)$ such that for any
$k\in \mathbb{T}_t$, the following FBS$\Delta$E admits a solution
$(X^{k,t,x},Z^{k,t,x})$
\begin{eqnarray}\label{system-adjoint}
\left\{\begin{array}{l}
X^{k,t,x}_{\ell+1}=\big{(}A_{k,\ell}X^{k,t,x}_\ell+\bar{A}_{k,\ell}
\mathbb{E}_kX^{k,t,x}_\ell+B_{k,\ell}u^{t,x,*}_\ell+\bar{B}_{k,\ell}
\mathbb{E}_ku^{t,x,*}_\ell+f_{k,\ell}\big{)}\\[1mm]
\hphantom{X^{k,t,x}_{\ell+1}=}+\big{(}C_{k,\ell}X^{k,t,x}_\ell+\bar{C}_{k,\ell}
\mathbb{E}_kX^{k,t,x}_\ell+D_{k,\ell}u^{t,x,*}_\ell+\bar{D}_{k,\ell}
\mathbb{E}_ku^{t,x,*}_\ell+d_{k,\ell}\big{)}w_\ell, \\[1mm]
Z_\ell^{k,t,x}=A_{k,\ell}^T\mathbb{E}_\ell Z_{\ell+1}^{k,t,x}+ \bar{A}_{k,\ell}^T
\mathbb{E}_kZ_{\ell+1}^{k,t,x}+C_{k,\ell}^T
\mathbb{E}_\ell(Z_{\ell+1}^{k,t,x}w_\ell)+\bar{C}_{k,\ell}^T
\mathbb{E}_k(Z_{\ell+1}^{k,t,x}w_\ell)\\[1mm]
\hphantom{Z_\ell^{k,t,x}=}+Q_{k,\ell}X_\ell^{k,t,x}+\bar{Q}_{k,\ell}
\mathbb{E}_kX_\ell^{k,t,x}+q_{k,\ell},\\[1mm]
X^{k,t,x}_k=X^{t,x,*}_k,~~\\[1mm]
Z_N^{k,t,x}=G_kX^{k,t,x}_N+\bar{G}_k
\mathbb{E}_kX^{k,t,x}_N+g_k,~~\ell\in \mathbb{T}_k
\end{array}
\right.
\end{eqnarray}
with the stationary condition
\begin{eqnarray}\label{stationary-condition}
0=(R_{k,k}+\bar{R}_{k,k})u_k^{t,x,*}+(B_{k,k}+\bar{B}_{k,k})^T
\mathbb{E}_kZ^{k,t,x}_{k+1}+(D_{k,k}+\bar{D}_{k,k})^T
\mathbb{E}_k(Z_{k+1}^{k,t,x}w_k)+\rho_{k,k},
\end{eqnarray}
and the convexity condition
\begin{eqnarray}\label{convex}
\inf_{\bar{u}_k\in L^2_\mathcal{F}(k;
\mathbb{R}^m)}\widehat{{J}}(k,0;\bar{u}_k) \geq 0.
\end{eqnarray}
Here, $\widehat{{J}}(k,0;\bar{u}_k)$ is given in (\ref{hat-J}), and
$X^{t,x,*}$ in (\ref{system-adjoint}) is given by
\begin{eqnarray}\label{open-loop-equilibrium-state}
\left\{
\begin{array}{l}
X^{t,x,*}_{k+1} = \big{[}(A_{k,k}+\bar{A}_{k,k})X^{t,x,*}_{k}+(B_{k,k}
+\bar{B}_{k,k})u^{t,x,*}_k+f_{k,k}\big{]}\\[1mm]
\hphantom{X^{t,x,*}_{k+1}=}+\big{[}(C_{k,k}+\bar{C}_{k,k})X^{t,x,*}_{k}
+(D_{k,k}+\bar{D}_{k,k})u^{t,x,*}_k+d_{k,k}\big{]}w_k,\\[1mm]
X^{t,x,*}_{t} = x,~~ k \in  \mathbb{T}_t.
\end{array}
\right.
\end{eqnarray}

Under any of above conditions, $(X^{t,x,*}, u^{t,x,*})$ given in (ii) is an open-loop equilibrium pair for the initial pair $(t,x)$.
\end{theorem}

\emph{Proof.}  See Appendix B. \endpf


To simplify the notations, let
\begin{eqnarray}\label{notation}
\left\{
\begin{array}{l}
\mathcal{A}_{k,\ell }=A_{k,\ell}+\bar{A}_{k,\ell},~~~
\mathcal{B}_{k,\ell }=B_{k,\ell}+\bar{B}_{k,\ell},~~~\\[1mm]
\mathcal{C}_{k,\ell }=C_{k,\ell}+\bar{C}_{k,\ell},~~~
\mathcal{D}_{k,\ell }=D_{k,\ell}+\bar{D}_{k,\ell},\\[1mm]
\mathcal{Q}_{k,\ell }=Q_{k,\ell}+\bar{Q}_{k,\ell},~~~
\mathcal{R}_{k,\ell }=R_{k,\ell}+\bar{R}_{k,\ell},~~~\\[1mm]
\mathcal{G}_k=G_k+\bar{G}_k,~~~k\in \mathbb{T}_t,~~
\ell\in \mathbb{T}_k.
\end{array}
\right.
\end{eqnarray}
Then, (\ref{open-loop-equilibrium-state}) can be rewritten as
\begin{eqnarray}\label{open-loop-equilibrium-state-0}
\left\{
\begin{array}{l}
X^{t,x,*}_{k+1} =\big{[}\mathcal{A}_{k,k}X^{t,x,*}_{k}
+\mathcal{B}_{k,k}u^{t,x,*}_k+f_{k,k}\big{]}+\big{[}
\mathcal{C}_{k,k}X^{t,x,*}_{k} +\mathcal{D}_{k,k}u^{t,x,*}_k+d_{k,k}\big{]}w_k,\\[1mm]
X^{t,x,*}_{t} = x,~~ k \in  \mathbb{T}_t.
\end{array}
\right.
\end{eqnarray}
For any $k\in \mathbb{T}_t$, from (\ref{system-adjoint}) it follows
that
\begin{eqnarray*}\label{system-adjoint-N-1}
\left\{\begin{array}{l}
X^{k,t,x}_{k+1}=\big{[}\mathcal{A}_{k,k}X^{k,t,x}_{k}
+\mathcal{B}_{k,k}u^{t,x,*}_{k}+f_{k,k}\big{]}+\big{[}
\mathcal{C}_{k,k}X^{k,t,x}_{k}+\mathcal{D}_{k,k}u^{t,x,*}_{k}+d_{k,k}\big{]}w_{k}, \\[2mm]
X^{k,t,x}_{k}=X^{t,x,*}_{k}.
\end{array}
\right.
\end{eqnarray*}
Therefore, we have %
\begin{eqnarray}\label{X}
X^{k,t,x}_{k+1}=X^{t,x,*}_{k+1},~~~k\in \mathbb{T}_t.
\end{eqnarray}

We now study the condition (\ref{convex}). The following result gives an
expression of $\widehat{{J}}(k,0;\bar{u}_k)$.

\begin{lemma}\label{Lemma-convex}
$\widehat{{J}}(k,0;\bar{u}_k)$ can be expressed as
\begin{eqnarray}\label{convex-2}
\widehat{{J}}(k,0;\bar{u}_k)=\bar{u}_k^T\big{(}\mathcal{R}_{k,k}
+\mathcal{B}_{k,k}^T\mathcal{P}_{k,k+1}\mathcal{B}_{k,k}+\mathcal{D}_{k,k}^TP_{k,k+1}\mathcal{D}_{k,k}\big{)}\bar{u}_k
\end{eqnarray}
with $P_{k,k+1}$ and $\mathcal{P}_{k,k+1}$ computed via
\begin{eqnarray}\label{P-0}
\left\{
\begin{array}{l}
P_{k,\ell}=Q_{k,\ell}+A^T_{k,\ell}P_{k,\ell+1}A_{k,\ell}+C^T_{k,\ell}
P_{k,\ell+1}C_{k,\ell},\\[1mm]
{\mathcal{P}}_{k,\ell}=\mathcal{Q}_{k,\ell}+\mathcal{A}^T_{k,\ell}
\mathcal{P}_{k,\ell+1}{\mathcal{A}}_{k,\ell}+\mathcal{C}^T_{k,\ell}
P_{k,\ell+1}{\mathcal{C}}_{k,\ell},\\[1mm]
P_{k,N}=G_k,~~{\mathcal{P}}_{k,N}={\mathcal{G}}_k, ~~\ell\in \mathbb{T}_k.
\end{array}
\right.
\end{eqnarray}
\end{lemma}

\emph{Proof.} From (\ref{system-y-k}), it follows that
\begin{eqnarray*}\label{system-Ey-k}
\left\{
\begin{array}{l}
\mathbb{E}_k{Y}^{k,\bar{u}_k}_{\ell+1}=\mathcal{A}_{k,\ell}
\mathbb{E}_k{Y}^{k,\bar{u}_k}_\ell,~~~\ell\in \mathbb{T}_{k+1},\\[1mm]
\mathbb{E}_kY^{k,\bar{u}_k}_{k+1}=\mathcal{B}_{k,k}\mathbb{E}_k\bar{u}_k,\\
\mathbb{E}_k{Y}^{k,\bar{u}_k}_k=0.
\end{array}
\right.
\end{eqnarray*}
By adding to and subtracting
\begin{eqnarray*}
&&\sum_{\ell=k}^{N-1}\mathbb{E}_k\Big{[}(Y^{k,\bar{u}_k}_{\ell+1})^T
P_{\ell+1}Y^{k,\bar{u}_k}_{\ell+1}-
(Y^{k,\bar{u}_k}_{\ell})^TP_{\ell}Y^{k,\bar{u}_k}_{\ell}\\
&&+(\mathbb{E}_kY^{k,\bar{u}_k}_{\ell+1})^T\bar{P}_{\ell+1}
\mathbb{E}_kY^{k,\bar{u}_k}_{\ell+1}-(\mathbb{E}_kY^{k,\bar{u}_k}_{\ell})^T
\bar{P}_{\ell}\mathbb{E}_kY^{k,\bar{u}_k}_{\ell}
\Big{]}
\end{eqnarray*}
from (\ref{hat-J}), we have
\begin{eqnarray}\label{convex-5}
&&\hspace{-2.2em}\widehat{{J}}(k,0;\bar{u}_k)=\sum_{\ell=k}^{N-1}\mathbb{E}_k\Big{\{}(Y_{\ell}^{k,\bar{u}_k})^TQ_{k,\ell}
Y_\ell^{k,\bar{u}_k}+(\mathbb{E}_kY_{\ell}^{k,\bar{u}_k})^T\bar{Q}_{k,\ell}
\mathbb{E}_kY_\ell^{k,\bar{u}_k}+ (Y^{k,\bar{u}_k}_{\ell+1})^TP_{\ell+1}
Y^{k,\bar{u}_k}_{\ell+1}\nonumber \\
&&\hspace{-2.2em}\hphantom{\widehat{{J}}(k,0;\bar{u}_k)=} - (Y^{k,\bar{u}_k}_{\ell})^TP_{\ell}Y^{k,
\bar{u}_k}_{\ell}+(\mathbb{E}_kY^{k,\bar{u}_k}_{\ell+1})^T
\bar{P}_{\ell+1}\mathbb{E}_kY^{k,\bar{u}_k}_{\ell+1}-(\mathbb{E}_kY^{k,
\bar{u}_k}_{\ell})^T\bar{P}_{\ell}\mathbb{E}_kY^{k,\bar{u}_k}_{\ell}\Big{\}}+\bar{u}_k^T\mathcal{R}_{k,k}\bar{u}_k\nonumber \\[1mm]
&&\hspace{-2.2em}\hphantom{\widehat{{J}}(k,0;\bar{u}_k)}=\sum_{\ell=k+1}^{N-1}\mathbb{E}_k\Big{\{}(\mathbb{E}_kY^{k,\bar{u}_k}_{\ell})^T
\big{(}\mathcal{Q}_{k,\ell}+\mathcal{A}^T_{k,\ell}\mathcal{P}_{k,\ell+1}
\mathcal{A}_{k,\ell}+\mathcal{C}^T_{k,\ell}{P}_{k,\ell+1}\mathcal{C}_{k,\ell}
-\mathcal{P}_{k,\ell} \big{)}\mathbb{E}_kY^{k,\bar{u}_k}_{\ell}\nonumber \\[1mm]
&&\hspace{-2.2em}\hphantom{\widehat{{J}}(k,0;\bar{u}_k)=}+(Y^{k,\bar{u}_k}_\ell
-\mathbb{E}_kY^{k,\bar{u}_k}_\ell)^T\big{(} Q_{k,\ell}+A_{k,\ell}^TP_{k,\ell+1}
A_{k,\ell}+C^T_{k,\ell}P_{k,\ell+1}C_{k,\ell}-P_{k,\ell} \big{)}(Y^{k,
\bar{u}_k}_\ell-\mathbb{E}_kY^{k,\bar{u}_k}_\ell)\Big{\}}\nonumber \\[1mm]
&&\hspace{-2.2em}\hphantom{\widehat{{J}}(k,0;\bar{u}_k)=}+\bar{u}_k^T\big{(}\mathcal{R}_{k,k}+\mathcal{B}_{k,k}^T
\mathcal{P}_{k,k+1}\mathcal{B}_{k,k}+\mathcal{D}_{k,k}^TP_{k,k+1}
\mathcal{D}_{k,k}\big{)}\bar{u}_k\nonumber \\[1mm]
&&\hspace{-2.2em}\hphantom{\widehat{{J}}(k,0;\bar{u}_k)}=\bar{u}_k^T\big{(}\mathcal{R}_{k,k}+\mathcal{B}_{k,k}^T
\mathcal{P}_{k,k+1}\mathcal{B}_{k,k}+\mathcal{D}_{k,k}^TP_{k,k+1}
\mathcal{D}_{k,k}\big{)}\bar{u}_k. \end{eqnarray}
This completes the proof. \endpf

Letting $u_k=0, \lambda=1$ in (\ref{appendix-A-J-0}), from Lemma
\ref{Lemma-convex} we have
\begin{eqnarray}\label{J-1}
&&\hspace{-4em}{J}(k,\zeta; (\bar{u}_k,u|_{\mathbb{T}_{k+1}}))=
\bar{u}_k^T\big{(}\mathcal{R}_{k,k}+\mathcal{B}_{k,k}^T\mathcal{P}_{k,k+1}
\mathcal{B}_{k,k}+\mathcal{D}_{k,k}^TP_{k,k+1}\mathcal{D}_{k,k}\big{)}
\bar{u}_k\nonumber\\[1mm]
&&\hspace{-4em}\hphantom{{J}(k,\zeta; (\bar{u}_k,u|_{\mathbb{T}_{k+1}}))=}+2
\Big{[}\rho_{k,k}+\mathcal{B}_{k,k}^T\mathbb{E}_kZ^{k,0}_{k+1}+\mathcal{D}_{k,k}^T
\mathbb{E}_k(Z_{k+1}^{k,0}w_k)\Big{]}^T\bar{u}_k+{J}(k,\zeta; (0,u|_{\mathbb{T}_{k+1}}))\nonumber\\[1mm]
&&\hspace{-4em}
\hphantom{{J}(k,\zeta; (\bar{u}_k,u|_{\mathbb{T}_{k+1}}))}\triangleq \langle M_{k,2} \bar{u}_k, \bar{u}_k\rangle+2\langle
M_{k,1}, \bar{u}_k \rangle + {J}(k,\zeta;
(0,u|_{\mathbb{T}_{k+1}})),
\end{eqnarray}
where $\langle \cdot , \cdot \rangle$ is the inner product on $\mathbb{R}^m$, and
$M_{k,2}$, $M_{k,1}$ are defined as
\begin{eqnarray}\label{M1-M2}
\left\{
\begin{array}{l}
M_{k,2}=\mathcal{R}_{k,k}+\mathcal{B}_{k,k}^T\mathcal{P}_{k,k+1}
\mathcal{B}_{k,k}+\mathcal{D}_{k,k}^TP_{k,k+1}\mathcal{D}_{k,k},\\[1mm]
M_{k,1}=\mathcal{B}_{k,k}^T\mathbb{E}_kZ^{k,0}_{k+1}+\rho_{k,k}+\mathcal{D}_{k,k}^T\mathbb{E}_k(Z_{k+1}^{k,0}w_k).
\end{array}
\right.
\end{eqnarray}
In (\ref{M1-M2}), $Z_{k+1}^{k,0}$ is computed via a version of (\ref{X-Z-1}) with $u_k$
replaced by $0$. It can be seen that $Z^{k,0}_{k+1}$ is a functional
of $\zeta$ and $u|_{\mathbb{T}_{k+1}}$.

Fixing $\zeta$ and $u|_{\mathbb{T}_{k+1}}$, ${J}(k,\zeta;
(\bar{u}_k,u|_{\mathbb{T}_{k+1}}))$ is a quadratic functional of $\bar{u}_k$, and is convex with respect to $\bar{u}_k$ if
$M_{k,2}\geq 0$. Throughout this paper, Problem (LQ) will be called
point-wisely convex at $(k,\zeta)$ (with $\zeta\in
L^2_{\mathcal{F}}(k; \mathbb{R}^n)$) if for fixed
$u|_{\mathbb{T}_{k+1}}$, $J(k, \zeta;
(\bar{u}_k,u|_{\mathbb{T}_{k+1}}))$ is convex with respect to $\bar{u}_k$. 
%
By this, Lemma \ref{Lemma-convex} and Theorem
\ref{Theorem-Equivalentce-open-loop}, we have the following result,
whose proof is omitted here.

\begin{proposition}\label{proposition-convex}
The following statements are equivalent.

(i) The convexity condition (\ref{convex}) is satisfied.

(ii) The following inequality holds
\begin{eqnarray}\label{convex-6}
M_{k,2}=\mathcal{R}_{k,k}+\mathcal{B}_{k,k}^T\mathcal{P}_{k,k+1}\mathcal{B}_{k,k}
+\mathcal{D}_{k,k}^TP_{k,k+1}\mathcal{D}_{k,k}\geq 0,
\end{eqnarray}
where $P_{k,k+1}$ and $\mathcal{P}_{k,k+1}$  are computed via (\ref{P-0}).

(iii) Problem (LQ) is point-wisely convex at $(k,\zeta)$ with some
$\zeta\in  L^2_{\mathcal{F}}(k; \mathbb{R}^n)$.

(iv) Problem (LQ) is point-wisely convex at $(k,\zeta)$ with any
$\zeta\in  L^2_{\mathcal{F}}(k; \mathbb{R}^n)$.

Furthermore, if Problem (LQ) admits an open-loop equilibrium pair
for the initial pair $(t,x)$, then for any $k\in \mathbb{T}_t$ and any
$\zeta \in L^2_{\mathcal{F}}(k; \mathbb{R}^n)$, Problem (LQ) is
point-wisely convex at $(k,\zeta)$.
\end{proposition}

Now let us switch to the stationary condition (\ref{stationary-condition}). The following lemma gives an expression of the backward state
$Z^{k,t,x}$.

\begin{lemma}\label{Lemma-Z}
Let $u^{t,x,*}_\ell=\Psi_\ell X^{t,x,*}_\ell +\alpha_\ell, \ell\in \mathbb{T}_{k}$, in (\ref{system-adjoint}) with $\Psi_\ell, \alpha_\ell, \ell\in \mathbb{T}_{k}$, being deterministic matrices. Then, the backward state $Z^{k,t,x}$ of (\ref{system-adjoint}) has the following expression
\begin{eqnarray}\label{Z}
&&\hspace{-2em}Z^{k,t,x}_\ell=P_{k,\ell}X^{k,t,x}_\ell+\bar{P}_{k,\ell}
\mathbb{E}_kX^{k,t,x}_\ell+T_{k,\ell}X^{t,x,*}_\ell+\bar{T}_{k,\ell}\mathbb{E}_kX^{t,x,*}_\ell+\pi_{k,\ell},~~~~\ell\in
\mathbb{T}_{k}.
\end{eqnarray}
Here, $\bar{P}_{k,\ell}=\mathcal{P}_{k,\ell}-P_{k,\ell}$ with $P_{k,\ell}, \mathcal{P}_{k,\ell}$ being computed via (\ref{P-0}),
and  $T_{k,\ell}, \bar{T}_{k,\ell}, \pi_{k,\ell}$ are given by
\begin{eqnarray}\label{T-lemma}
\left\{
\begin{array}{l}
T_{k,\ell}=A^T_{k,\ell}T_{k,\ell+1}\mathcal{A}_{\ell,\ell}
+C^T_{k,\ell}T_{k,\ell+1}\mathcal{C}_{\ell,\ell}\\[1mm]
\hphantom{T_{k,\ell}=}+\Big{(}A^T_{k,\ell}P_{k,\ell+1}B_{k,\ell}
+A^T_{k,\ell}T_{k,\ell+1}\mathcal{B}_{\ell,\ell}
+C^T_{k,\ell}P_{k,\ell+1}D_{k,\ell}+C^T_{k,\ell}T_{k,\ell+1}
\mathcal{D}_{\ell,\ell}\Big{)}\Psi_{\ell},\\[2mm]
\bar{T}_{k,\ell}=A^T_{k,\ell}\bar{T}_{k,\ell+1}\mathcal{A}_{\ell,\ell}
+\bar{A}^T_{k,\ell}\mathcal{T}_{k,\ell+1}\mathcal{A}_{\ell,\ell}+\bar{C}^T_{k,\ell}T_{k,\ell+1}\mathcal{C}_{\ell,\ell}
\\[1mm]
\hphantom{\bar{T}_{k,\ell}=}+\Big{(}A^T_{k,\ell}P_{k,\ell+1}
\bar{B}_{k,\ell}+A^T_{k,\ell}\bar{P}_{k,\ell+1}
\mathcal{{B}}_{k,\ell}+A^T_{k,\ell}\bar{T}_{k,\ell+1}
\mathcal{B}_{\ell,\ell}+C^T_{k,\ell}P_{k,\ell+1}\bar{D}_{k,\ell}\\[1mm]
\hphantom{\bar{T}_{k,\ell}=}+\bar{A}^T_{k,\ell}\mathcal{P}_{k,\ell+1}
\mathcal{B}_{k,\ell}+\bar{A}^T_{k,\ell}\mathcal{T}_{k,\ell+1}
\mathcal{B}_{\ell,\ell}+\bar{C}^T_{k,\ell}P_{k,\ell+1}
\mathcal{D}_{k,\ell}+\bar{C}^T_{k,\ell}T_{k,\ell+1}
\mathcal{D}_{\ell,\ell}\Big{)}\Psi_{\ell},\\[1mm]
T_{k,N}=0,~~\bar{T}_{k,N}=0,\\[1mm]
\ell\in \mathbb{T}_{k},
\end{array}
\right.
\end{eqnarray}
and
\begin{eqnarray}\label{pi-lemma}
\left\{
\begin{array}{l}
\pi_{k,\ell}=\mathcal{A}^T_{k,\ell}\mathcal{P}_{k,\ell+1}\big{(}\mathcal{B}_{k,\ell}
\alpha_{\ell}+f_{k,\ell}\big{)}+\mathcal{A}^T_{k,\ell}\mathcal{T}_{k,\ell+1}
\big{(}\mathcal{B}_{\ell,\ell} \alpha_{\ell}+f_{\ell,\ell}\big{)}+\mathcal{A}^T_{k,\ell}\pi_{k,\ell+1}\\[1mm]
\hphantom{\pi_{k,\ell}=}+\mathcal{C}^T_{k,\ell}P_{k,\ell+1}\big{(}
\mathcal{D}_{k,\ell}\alpha_{\ell}+d_{k,\ell}\big{)}+\mathcal{C}^T_{k,\ell}
T_{k,\ell+1}\big{(}\mathcal{D}_{\ell,\ell}\alpha_{\ell} +d_{\ell,\ell}\big{)}+q_{k,\ell},  \\[1mm]
\pi_{k,N}=g_k,~~\\[1mm]
\ell\in \mathbb{T}_{k}
\end{array}
\right.
\end{eqnarray}
with $\mathcal{T}_{k,\ell}=T_{k,\ell}+\bar{T}_{k,\ell}, \ell\in \mathbb{T}_k$.

\end{lemma}

\emph{Proof}. See Appendix C. \endpf

Noting that $P_{k,\ell}$ and
$\mathcal{P}_{k,\ell}$ are symmetric, $T_{k,\ell}$ and $\bar{T}_{k,\ell}$
are generally nonsymmetric as $\mathcal{A}_{\ell,\ell}, \mathcal{B}_{\ell,\ell},
\mathcal{C}_{\ell,\ell}$ and $\mathcal{D}_{\ell,\ell}$ appear in the
expressions of $T_{k,\ell}$ and $\bar{T}_{k,\ell}$.
Recall the pseudo-inverse of a matrix. By \cite{Penrose}, for a
given matrix $M\in \mathbb{R}^{n\times m}$, there exists a unique
matrix in $\mathbb{R}^{m\times n}$ denoted by $M^\dagger$ such that
\begin{eqnarray}
\left\{
\begin{array}{l}
MM^\dagger M=M,~~ M^\dagger M M^\dagger=M^\dagger,\\
(MM^\dagger)^T=MM^\dagger, ~~(M^\dagger M)^T=M^\dagger M.
\end{array}
\right.
\end{eqnarray}
This $M^\dagger$ is called the Moore-Penrose inverse  of $M$. The
following lemma is from \cite{Ait-Chen-Zhou-2002}.

\begin{lemma}\label{Lemma-matrix-equation}
Let matrices $L$, $M$ and $N$ be given with appropriate size. Then,
$LXM=N$ has a solution $X$ if and only if $LL^\dagger NMM^\dagger=N$.
Moreover, the solution of $LXM=N$ can be expressed as
$X=L^\dagger NM^\dagger+Y-L^\dagger LYMM^\dagger$,
where $Y$ is a matrix with appropriate size.
\end{lemma}

Based on above results, we have the following theorem.

\begin{theorem}\label{Theorem-Necessary}
Given $t\in \mathbb{T}$ and $x\in L^2_\mathcal{F}(t; \mathbb{R}^n)$, the
following statements are equivalent.

(i) There exists an open-loop equilibrium pair of Problem (LQ) for the initial pair $(t,x)$.

(ii) The set of LDEs
\begin{eqnarray}\label{P}
\left\{
\begin{array}{l}
\left\{
\begin{array}{l}
P_{k,\ell}=Q_{k,\ell}+A^T_{k,\ell}P_{k,\ell+1}A_{k,\ell}+C^T_{k,\ell}P_{k,\ell+1}C_{k,\ell},\\[1mm]
{\mathcal{P}}_{k,\ell}=\mathcal{Q}_{k,\ell}+\mathcal{A}^T_{k,\ell}\mathcal{P}_{k,\ell+1}{\mathcal{A}}_{k,\ell}+\mathcal{C}^T_{k,\ell}P_{k,\ell+1}{\mathcal{C}}_{k,\ell},\\[1mm]
P_{k,N}=G_k,~~\mathcal{{P}}_{k,N}=\mathcal{{G}}_k,~~\ell\in \mathbb{T}_k,
\end{array}
\right. \\[1mm]
\mathcal{R}_{k,k}+\mathcal{B}_{k,k}^T\mathcal{P}_{k,k+1}\mathcal{B}_{k,k}
+\mathcal{D}_{k,k}^TP_{k,k+1}\mathcal{D}_{k,k}\geq 0,\\[1mm]
k\in \mathbb{T}_t
\end{array}
\right.
\end{eqnarray}
is solvable in the sense
\begin{eqnarray}\label{convex-4}
\mathcal{R}_{k,k}+\mathcal{B}_{k,k}^T\mathcal{P}_{k,k+1}\mathcal{B}_{k,k}
+\mathcal{D}_{k,k}^TP_{k,k+1}\mathcal{D}_{k,k}\geq 0, k\in
\mathbb{T}_t,
\end{eqnarray}
and the following condition
\begin{eqnarray}\label{Theorem-Necessary-equality}
\big{(}I-{\mathcal{W}}_k \mathcal{W}_k^\dagger\big{)}\big{(}\mathcal{H}_kX^{t,x,*}_k+\beta_k\big{)}=0,~~k\in \mathbb{T}_t
\end{eqnarray}
is satisfied. Here,
\begin{eqnarray}\label{Theorem-Necessary-open-loop-equilibrium-state}
\left\{
\begin{array}{l}
X^{t,x,*}_{k+1} = \big{(}\mathcal{A}_{k,k}-\mathcal{B}_{k,k}\mathcal{W}_{k}^\dagger \mathcal{H}_{k}\big{)}X^{t,x,*}_{k}-\mathcal{B}_{k,k}\mathcal{W}_{k}^\dagger \beta_{k}+f_{k,k}\\[1mm]
\hphantom{X^{t,x,*}_{k+1}=}+\big{[}\big{(}\mathcal{C}_{k,k}-\mathcal{D}_{k,k}\mathcal{W}_{k}^\dagger \mathcal{H}_{k}\big{)}X^{t,x,*}_{k}-\mathcal{D}_{k,k}\mathcal{W}_{k}^\dagger \beta_{k}+d_{k,k}\big{]}w_k,\\[1mm]
X^{t,x,*}_{t} = x,~~ k \in  \mathbb{T}_t,
\end{array}
\right.
\end{eqnarray}
and
\begin{eqnarray}\label{feedback-gain-k}
\left\{
\begin{array}{l}
\mathcal{W}_{k}=\mathcal{R}_{k,k}+\mathcal{B}_{k,k}^T\big{(}\mathcal{P}_{k,k+1}
+\mathcal{T}_{k,k+1}\big{)}\mathcal{B}_{k,k}+\mathcal{D}_{k,k}^T\big{(}{P}_{k,k+1}
+T_{k,k+1}\big{)}\mathcal{D}_{k,k},\\[1mm]
\mathcal{H}_{k}=\mathcal{B}_{k,k}^T\big{(}\mathcal{P}_{k,k+1}+\mathcal{T}_{k,k+1}\big{)}
\mathcal{A}_{k,k}+\mathcal{D}_{k,k}^T\big{(}{P}_{k,k+1}+T_{k,k+1}\big{)}\mathcal{C}_{k,k},\\[1mm]
\beta_{k}=\mathcal{B}_{k,k}^T\big{[}\big{(}\mathcal{{P}}_{k,k+1}+\mathcal{T}_{k,k+1}
\big{)}f_{k,k}+\pi_{k,k+1}\big{]}+\mathcal{D}_{k,k}^T\big{(}{P}_{k,k+1}+T_{k,k+1}
\big{)}d_{k,k}+\rho_{k,k},\\[1mm] k\in \mathbb{T}_t
\end{array}
\right.
\end{eqnarray}
with
\begin{eqnarray}\label{T}
\left\{
\begin{array}{l}
\left\{
\begin{array}{l}
T_{k,\ell}=A^T_{k,\ell}T_{k,\ell+1}\mathcal{A}_{\ell,\ell}
+C^T_{k,\ell}T_{k,\ell+1}\mathcal{C}_{\ell,\ell}\\[1mm]
\hphantom{T_{k,\ell}=}-\Big{(}A^T_{k,\ell}P_{k,\ell+1}B_{k,\ell}
+A^T_{k,\ell}T_{k,\ell+1}\mathcal{B}_{\ell,\ell}+C^T_{k,\ell}P_{k,\ell+1}D_{k,\ell}+C^T_{k,\ell}T_{k,\ell+1}
\mathcal{D}_{\ell,\ell}\Big{)}\mathcal{W}_{\ell}^\dagger \mathcal{H}_{\ell},\\[1mm]
{\mathcal{T}}_{k,\ell}=\mathcal{A}^T_{k,\ell}{\mathcal{T}}_{k,\ell+1}
\mathcal{A}_{\ell,\ell}+{\mathcal{C}}^T_{k,\ell}T_{k,\ell+1}
\mathcal{C}_{\ell,\ell}\\[1mm]
\hphantom{\bar{T}_{k,\ell}=}-\Big{(}\mathcal{A}^T_{k,\ell}
\mathcal{P}_{k,\ell+1}{\mathcal{B}}_{k,\ell}
+\mathcal{A}^T_{k,\ell}{\mathcal{T}}_{k,\ell+1}\mathcal{B}_{\ell,\ell}+\mathcal{C}^T_{k,\ell}P_{k,\ell+1}{\mathcal{D}}_{k,\ell}+{\mathcal{C}}^T_{k,\ell}
T_{k,\ell+1}\mathcal{D}_{\ell,\ell}\Big{)}\mathcal{W}_{\ell}^\dagger \mathcal{H}_{\ell},\\[1mm]
T_{k,N}=0,~~{\mathcal{T}}_{k,N}=0, \\[1mm]
\ell\in \mathbb{T}_{k},
\end{array}
\right.\\[1mm]
k\in \mathbb{T}_t,
\end{array}
\right.
\end{eqnarray}
and
\begin{eqnarray}\label{pi}
\left\{
\begin{array}{l}
\left\{
\begin{array}{l}
\pi_{k,\ell}=\mathcal{A}^T_{k,\ell}\mathcal{P}_{k,\ell+1}\big{(}f_{k,\ell}
-\mathcal{B}_{k,\ell}\mathcal{W}_{\ell}^\dagger \beta_{\ell}\big{)}+\mathcal{A}^T_{k,\ell}\mathcal{T}_{k,\ell+1}\big{(}f_{\ell,\ell}
-\mathcal{B}_{\ell,\ell} \mathcal{W}_{\ell}^\dagger \beta_{\ell}\big{)}
\\[1mm]
\hphantom{\pi_{k,\ell}=}+\mathcal{C}^T_{k,\ell}P_{k,\ell+1}\big{(}d_{k,\ell}
-\mathcal{D}_{k,\ell}\mathcal{W}_{\ell}^\dagger \beta_{\ell}\big{)}+\mathcal{C}^T_{k,\ell}T_{k,\ell+1}\big{(}d_{\ell,\ell}-\mathcal{D}_{\ell,\ell}
\mathcal{W}_{\ell}^\dagger \beta_{\ell}\big{)}\\[1mm]
\hphantom{\pi_{k,\ell}=}+\mathcal{A}^T_{k,\ell}\pi_{k,\ell+1}+q_{k,\ell},\\[1mm]
\pi_{k,N}=g_k,
\end{array}
\right.\\[1mm]
k\in \mathbb{T}_t.
\end{array}
\right.
\end{eqnarray}
Furthermore, we have
\begin{eqnarray}\label{Z-2}
&&\hspace{-2em}Z^{k,t,x}_\ell=P_{k,\ell}\big{(}X^{k,t,x}_\ell-\mathbb{E}_kX^{k,t,x}_\ell
\big{)}+\mathcal{{P}}_{k,\ell}\mathbb{E}_kX^{k,t,x}_\ell+T_{k,\ell}\big{(}X^{t,x,*}_\ell-\mathbb{E}_k
X^{t,x,*}_\ell \big{)}\nonumber \\[1mm]
&&\hspace{-2em}\hphantom{Z^{k,t,x}_\ell=}+{\mathcal{T}}_{k,\ell}\mathbb{E}_kX^{t,x,*}_\ell
+\pi_{k,\ell},~~~~\ell\in \mathbb{T}_{k}.
\end{eqnarray}

Under any of above conditions, an open-loop equilibrium control is given by
\begin{eqnarray}\label{open-loop-equilibrium-k}
u_{k}^{t,x,*}=-\mathcal{W}_{k}^\dagger \mathcal{H}_{k} X_{k}^{t,x,*}
-\mathcal{W}_{k}^\dagger \beta_{k},~~k\in \mathbb{T}_t
\end{eqnarray}
with $X^{t,x,*}$ given by (\ref{Theorem-Necessary-open-loop-equilibrium-state}).

\end{theorem}

\emph{Proof}. (i)$\Rightarrow$(ii). Let $u^{t,x,*}$ be an
open-loop equilibrium control. 
From Theorem \ref{Theorem-Equivalentce-open-loop} and Proposition \ref{proposition-convex}, we have the solvability of (\ref{P}). Furthermore,
for any $k\in \mathbb{T}_t$, the FBS$\Delta$E (\ref{system-adjoint}) admits a  solution, and (\ref{stationary-condition}) holds.
As %
\begin{eqnarray*}
Z_{N}^{N-1,t,x}=G_{N-1}X^{N-1,t,x}_{N}+\bar{G}_{N-1}\mathbb{E}_{N-1}X^{N-1,t,x}_{N}+g_{N-1},
\end{eqnarray*}
from (\ref{stationary-condition}) and (\ref{X}) we have
\begin{eqnarray*}
&&0=\mathcal{R}_{N-1,N-1}u_{N-1}^{t,x,*}+\mathcal{B}_{N-1,N-1}^T
\mathcal{G}_{N-1}\mathbb{E}_{N-1}X_{N}^{t,x,*}\\[1mm]
&&\hphantom{0=}+\mathcal{D}_{N-1,N-1}^T{G}_{N-1}\mathbb{E}_{N-1}(X_{N}^{t,x,*}w_{N-1})+\mathcal{B}_{N-1,N-1}^Tg_{N-1}+\rho_{N-1,N-1}.
\end{eqnarray*}
Note that $X^{t,x,*}$ is given in (\ref{open-loop-equilibrium-state}).
Then, substituting $X^{t,x,*}_{N-1}$ into the above equation, from  Lemma \ref{Lemma-matrix-equation} we have
\begin{eqnarray}\label{open-loop-equilibrium-N-1--00}
&&\hspace{-2em}u_{N-1}^{t,x,*}=-\mathcal{W}_{N-1}^\dagger \mathcal{H}_{N-1} X_{N-1}^{t,x,*}
-\mathcal{W}_{N-1}^\dagger \beta_{N-1}\nonumber \\[1mm]
&&\hspace{-2em}\hphantom{u_{N-1}^{t,x,*}}\triangleq \Psi_{N-1}X_{N-1}^{t,x,*}+\alpha_{N-1},
\end{eqnarray}
and
\begin{eqnarray*}
\big{(}I-{\mathcal{W}}_{N-1} \mathcal{W}_{N-1}^\dagger\big{)}\big{(}\mathcal{H}_{N-1}X^{t,x,*}_{N-1}+\beta_{N-1}\big{)}=0,
\end{eqnarray*}
where
\begin{eqnarray*}\label{feedback-gain-N-1}
\left\{
\begin{array}{l}
\mathcal{W}_{N-1}=\mathcal{R}_{N-1,N-1}+\mathcal{B}_{N-1,N-1}^T\mathcal{G}_{N-1}
\mathcal{B}_{N-1,N-1}+\mathcal{D}_{N-1,N-1}^T{G}_{N-1}\mathcal{D}_{N-1,N-1},\\[1mm]
\mathcal{H}_{N-1}=\mathcal{B}_{N-1,N-1}^T\mathcal{G}_{N-1}\mathcal{A}_{N-1,N-1}+\mathcal{D}_{N-1,N-1}^T{G}_{N-1}\mathcal{C}_{N-1,N-1},\\[1mm]
\beta_{N-1}=\mathcal{B}_{N-1,N-1}^T\big{[}\mathcal{G}_{N-1}f_{N-1,N-1}+g_{N-1}
\big{]}+\mathcal{D}_{N-1,N-1}^T{G}_{N-1}d_{N-1,N-1}+\rho_{N-1,N-1}.
\end{array}
\right.
\end{eqnarray*}
Noting (\ref{X}) and Lemma \ref{Lemma-Z}, we have
\begin{eqnarray*}
&&Z^{N-2,t,x}_{N-1}=\big{(}P_{N-2,N-1}+T_{N-2,N-1}\big{)}X^{t,x,*}_{N-1}+\big{(}\bar{P}_{N-2,N-1}+\bar{T}_{N-2,N-1}\big{)}\mathbb{E}_{N-2}
X^{t,x,*}_{N-1}+\pi_{N-2,N-1},
\end{eqnarray*}
where $T_{N-2,N-1}, \bar{T}_{N-2,N-1}$ are computed via (\ref{T-lemma}) with $\Psi_{N-1}$ and $\alpha_{N-1}$ being given in (\ref{open-loop-equilibrium-N-1--00}).
From (\ref{stationary-condition}), we have for $k=N-2$
\begin{eqnarray*}
&&0=\mathcal{R}_{N-2,N-2}u_{N-2}^{t,x,*}+\mathcal{B}_{N-2,N-2}^T
\Big{[}\big{(}\mathcal{P}_{N-2,N-1}+\mathcal{T}_{N-2,N-1}\big{)}
\mathbb{E}_{N-2}X_{N-1}^{t,x,*} +\pi_{N-2,N-1}\Big{]}\nonumber  \\[1mm]
&&\hphantom{0=}+\mathcal{D}_{N-2,N-2}^T\big{(}{P}_{N-2,N-1}
+T_{N-2,N-1}\big{)} \mathbb{E}_{N-2}\big{(}X_{N-1}^{t,x,*}w_{N-2}\big{)}+\rho_{N-2,N-2}.
\end{eqnarray*}
Substituting $X^{t,x,*}_{N-2}$ into above equation, by Lemma
\ref{Lemma-matrix-equation} we have
\begin{eqnarray*}\label{open-loop-equilibrium-N-1}
&&u_{N-2}^{t,x,*}=-\mathcal{W}_{N-2}^\dagger \mathcal{H}_{N-2} X_{N-2}^{t,x,*}
-\mathcal{W}_{N-2}^\dagger \beta_{N-2}\\[1mm]
&&\hphantom{u_{N-2}^{t,x,*}}\triangleq \Psi_{N-2}X_{N-2}^{t,x,*}+\alpha_{N-2},
\end{eqnarray*}
and
\begin{eqnarray*}
\big{(}I-{\mathcal{W}}_{N-2} \mathcal{W}_{N-2}^\dagger\big{)}\big{(}\mathcal{H}_{N-2}X^{t,x,*}_{N-2}+\beta_{N-2}\big{)}=0,
\end{eqnarray*}
where
\begin{eqnarray*}\label{feedback-gain-N-2}
\left\{\hspace{-1.5mm}
\begin{array}{l}
\mathcal{W}_{N-2}=\mathcal{R}_{N-2,N-2}+\mathcal{B}_{N-2,N-2}^T\big{(}
\mathcal{P}_{N-2,N-1}+\mathcal{T}_{N-2,N-1}\big{)}\mathcal{B}_{N-2,N-2}\\[1mm]
\hphantom{\mathcal{W}_{N-2}=}+\mathcal{D}_{N-2,N-2}^T \big{(}{P}_{N-2,N-1}
+T_{N-2,N-1}\big{)}\mathcal{D}_{N-2,N-2},\\[1mm]
\mathcal{H}_{N-2}=\mathcal{B}_{N-2,N-2}^T\big{(}\mathcal{P}_{N-2,N-1}+\mathcal{T}_{N-2,N-1}\big{)}\mathcal{A}_{N-2,N-2}\\[1mm]
\hphantom{\mathcal{H}_{N-2}=}
+\mathcal{D}_{N-2,N-2}^T\big{(}{P}_{N-2,N-1}+T_{N-2,N-1}\big{)} \mathcal{C}_{N-2,N-2},\\[1mm]
\beta_{N-2}=\mathcal{B}_{N-2,N-2}^T\big{[}\big{(}\mathcal{{P}}_{N-2,N-1}
+\mathcal{T}_{N-2,N-1}\big{)}f_{N-2,N-2}\\[1mm]
\hphantom{\beta_{N-2}=}+\mathcal{D}_{N-2,N-2}^T\big{(}{P}_{N-2,N-1}
+T_{N-2,N-1}\big{)}d_{N-2,N-2}\\[1mm]
\hphantom{\beta_{N-2}=}
+\pi_{N-2,N-1}\big{]}+\rho_{N-2,N-2}.
\end{array}
\right.
\end{eqnarray*}
Backwardly repeating above procedure, by Lemma \ref{Lemma-Z} we
can get (\ref{T}), (\ref{pi}) and (\ref{open-loop-equilibrium-k}).

(ii)$\Rightarrow$(i). By Proposition \ref{proposition-convex}, Lemma \ref{Lemma-matrix-equation} and reversing the proof of (i)$\Rightarrow$(ii), we can  complete the proof. \endpf

Now let Problem (LQ) for the initial pair $(t,x)$ admit an open-loop
equilibrium pair.
For a $\xi\neq x$ with
$\xi\in L^2_{\mathcal{F}}(t; \mathbb{R}^n)$, we can construct a
control of the form (\ref{open-loop-equilibrium-k}) as
\begin{eqnarray}\label{open-loop-equilibrium-k-8}
u_{k}^{t,\xi,*}=-\mathcal{W}_{k}^\dagger \mathcal{H}_{k} X_{k}^{t,\xi,*}
-\mathcal{W}_{k}^\dagger \beta_{k},~~k\in \mathbb{T}_t,
\end{eqnarray}
where
\begin{eqnarray*}
\left\{
\begin{array}{l}
{X}^{t,\xi,*}_{k+1} = \big{[}\big{(}\mathcal{A}_{k,k}-\mathcal{B}_{k,k}
\mathcal{W}^\dagger_k\mathcal{H}_k\big{)}{X}^{t,\xi,*}_{k}-\mathcal{B}_{k,k}\mathcal{W}^\dagger_k\beta_k+f_{k,k}\big{]}\\[1mm]
\hphantom{{X}^{t,\xi,*}_{k+1} =}+\big{[}\big{(}\mathcal{A}_{k,k} -\mathcal{B}_{k,k}\mathcal{W}^\dagger_k\mathcal{H}_k\big{)}{X}^{t,\xi,*}_{k}-\mathcal{D}_{k,k}\mathcal{W}^\dagger_k\beta_k+d_{k,k}\big{]}w_k,\\[1mm]
{X}^{t,\xi,*}_{t} = \xi,~~ k \in  \mathbb{T}_t,
\end{array}
\right.
\end{eqnarray*}
or equivalently,
\begin{eqnarray}\label{open-loop-equilibrium-state-2}
\left\{
\begin{array}{l}
X^{t,\xi,*}_{k+1} = \big{[}\mathcal{A}_{k,k}X^{t,\xi,*}_{k}+\mathcal{B}_{k,k}u^{t,\xi,*}_k+f_{k,k}\big{]}\\[1mm]
\hphantom{X^{t,\xi,*}_{k+1} = }+\big{[}\mathcal{C}_{k,k}X^{t,\xi,*}_{k} + \mathcal{D}_{k,k}u^{t,\xi,*}_k+d_{k,k}\big{]}w_k,\\[1mm]
X^{t,\xi,*}_{t} = \xi,~~ k \in  \mathbb{T}_t.
\end{array}
\right.
\end{eqnarray}
Though similarly defined as $(X^{t,{x},*}, u^{t,{x},*})$, we cannot
assert  that $(X^{t,\xi,*}, u^{t,\xi,*})$ is an open-loop
equilibrium pair of Problem (LQ) for the initial pair $(t,\xi)$. In fact, (\ref{Theorem-Necessary-equality}) reads as
\begin{eqnarray}\label{open-loop-equilibrium-W}
\mathcal{W}_{k}\mathcal{W}_{k}^\dagger\big{(}\mathcal{H}_{k} X_{k}^{t,x,*}
+ \beta_{k} \big{)}=\mathcal{H}_{k} X_{k}^{t,x,*}+ \beta_{k},~~k\in \mathbb{T}_t.
\end{eqnarray}
If $(X^{t,\xi,*}, u^{t,\xi,*})$ was an open-loop equilibrium pair for
the initial pair $(t,\xi)$, then there would be
\begin{eqnarray}\label{open-loop-equilibrium-W-2}
\mathcal{W}_{k}\mathcal{W}_{k}^\dagger\big{(}\mathcal{H}_{k} X_{k}^{t,\xi,*}
+ \beta_{k} \big{)}=\mathcal{H}_{k} X_{k}^{t,\xi,*}+ \beta_{k},~~k\in \mathbb{T}_t.
\end{eqnarray}
However, generally speaking, (\ref{open-loop-equilibrium-W-2})
cannot be deduced from (\ref{open-loop-equilibrium-W}). In fact, we
have
\begin{eqnarray*}
&&\hspace{-2em}\mathcal{W}_{k}\mathcal{W}_{k}^\dagger\big{(}\mathcal{H}_{k} X_{k}^{t,\xi,*}
+ \beta_{k} \big{)}=\mathcal{W}_{k}\mathcal{W}_{k}^\dagger\big{(}\mathcal{H}_{k}
X_{k}^{t,x,*}+ \beta_{k} \big{)}+\mathcal{W}_{k}\mathcal{W}_{k}^\dagger
\mathcal{H}_{k}\big{(} X_{k}^{t,\xi,*}- X_{k}^{t,x,*} \big{)}\\[1mm]
&&\hspace{-2em}\hphantom{\mathcal{W}_{k}\mathcal{W}_{k}^\dagger\big{(}\mathcal{H}_{k} X_{k}^{t,\xi,*}
+ \beta_{k} \big{)}}=\mathcal{H}_{k} X_{k}^{t,\xi,*}
+ \beta_{k}+\big{(}\mathcal{W}_{k}\mathcal{W}_{k}^\dagger\mathcal{H}_{k}
- \mathcal{H}_{k} \big{)}\big{(} X_{k}^{t,\xi,*}- X_{k}^{t,x,*} \big{)},
\end{eqnarray*}
which is different from $\mathcal{H}_{k} X_{k}^{t,\xi,*}+ \beta_{k}$
in general. Therefore, under the condition that Problem (LQ) has an open-loop equilibrium control for an initial pair $(t,x)$, we cannot assert the existence of open-loop equilibrium control for other initial pairs.

If (\ref{P}) is solvable, we have from Proposition
\ref{proposition-convex} that for any $(k,
\zeta)$ ($k\in \mathbb{T}_t$, $\zeta\in L^2_{\mathcal{F}}(k;
\mathbb{R}^n)$), ${J}(k,\zeta; (\bar{u}_k,u|_{\mathbb{T}_{k+1}}))$
is convex with respect to $\bar{u}_k$. By (\ref{J-1}), ${J}(k,\zeta;
(\bar{u}_k,u|_{\mathbb{T}_{k+1}}))$ can be rewritten as
\begin{eqnarray}\label{J-2}
&&\hspace{-2em}{J}(k,\zeta; (\bar{u}_k,u|_{\mathbb{T}_{k+1}}))= \langle M_{k,2}
\bar{u}_k, \bar{u}_k\rangle + {J}(k,\zeta; (0,u|_{\mathbb{T}_{k+1}}))+2\langle M_{k,1}(\zeta, u|_{\mathbb{T}_{k+1}}),
\bar{u}_k \rangle.
\end{eqnarray}
Here, we have used $M_{k,1}(\zeta, u|_{\mathbb{T}_{k+1}})$ instead
of $M_{k,1}$ to emphasize the dependence on $(\zeta,
u|_{\mathbb{T}_{k+1}})$. Only with the convexity condition, we cannot get the
existence of the minima of a quadratic functional like (\ref{J-2}).
To see more about this, let us consider a perturbation of the control weighting
matrices. Precisely, for $\varepsilon>0$ and $k\in \mathbb{T}_t$,
introduce the following cost functional
\begin{eqnarray}\label{J-3}
&&\hspace{-4em}{J}_\varepsilon(k,\zeta; (\bar{u}_k,u|_{\mathbb{T}_{k+1}}))={J}(k,\zeta; (\bar{u}_k,u|_{\mathbb{T}_{k+1}}))+\varepsilon \mathbb{E}\big{[}\bar{u}_k^T\bar{u}_k \big{]}\nonumber \\[1mm]
&&\hspace{-4em}\hphantom{{J}_\varepsilon(k,\zeta; (\bar{u}_k,u|_{\mathbb{T}_{k+1}}))}
= \langle \big{(}M_{k,2}+\varepsilon I\big{)} \bar{u}_k, \bar{u}_k\rangle+2\langle M_{k,1}(\zeta, u|_{\mathbb{T}_{k+1}}), \bar{u}_k \rangle \nonumber\\[1mm]
&&\hspace{-4em}\hphantom{{J}_\varepsilon(k,\zeta; (\bar{u}_k,u|_{\mathbb{T}_{k+1}}))=}+ {J}
(k,\zeta; (0,u|_{\mathbb{T}_{k+1}})).
\end{eqnarray}
Then, it holds that
\begin{eqnarray*}
&&\hspace{-2em}M_{k,2}^\varepsilon \triangleq M_{k,2}+\varepsilon I= \mathcal{R}_{k,k}
+\varepsilon I+\mathcal{B}_{k,k}^T\mathcal{P}_{k,k+1}\mathcal{B}_{k,k}
+\mathcal{D}_{k,k}^TP_{k,k+1}\mathcal{D}_{k,k}\geq \varepsilon I.
\end{eqnarray*}
By simple calculations, we have
\begin{eqnarray}\label{argmin}
&&\hspace{-2em}\bar{u}_k^*=arg \min_{\bar{u}_k\in L^2_{\mathcal{F}}(k; \mathbb{R}^m)}{J}_\varepsilon
(k,\zeta; (\bar{u}_k,u|_{\mathbb{T}_{k+1}}))\nonumber \\[1mm]
&&\hspace{-2em}\hphantom{\bar{u}_k^*}=-(M_{k,2}^\varepsilon)^{-1}
M_{k,1}(\zeta, u|_{\mathbb{T}_{k+1}}), ~~k\in \mathbb{T}_t.
\end{eqnarray}
In what follows, the version of Problem (LQ) corresponding to
$\{J_\varepsilon(k, \,\cdot \,, \,\cdot \,), k\in \mathbb{T}_t\}$ will be denoted as
Problem (LQ)$_\varepsilon$, for which we can adopt a
backward procedure to derive the open-loop equilibrium control.
Specifically, letting $k={N-1}$ in (\ref{J-3}) and by
(\ref{open-loop-equilibrium}) (\ref{argmin}), we have
\begin{eqnarray}\label{argmin-N-1}
{u}_{N-1}^{\varepsilon,t,x,*}=-(M_{N-1,2}^\varepsilon)^{-1}M_{N-1,1}
(X_{N-1}^{\varepsilon,t,x,*})
\end{eqnarray}
with the process $X^{\varepsilon,t,x,*}$ being determined below.
Substituting (\ref{argmin-N-1}) into (\ref{open-loop-equilibrium}),
from (\ref{argmin}) we have
\begin{eqnarray*}\label{argmin-N-2}
{u}_{N-2}^{\varepsilon,t,x,*}=-(M_{N-2,2}^\varepsilon)^{-1}M_{N-2,1}
(X_{N-2}^{\varepsilon,t,x,*}, u^{\varepsilon,t,x,*}_{N-1}).
\end{eqnarray*}
Repeating above procedure backwardly, we get the following
open-loop equilibrium pair $(X^{\varepsilon,t,x,*},
u^{\varepsilon,t,x,*})$:
\begin{eqnarray*}
{u}_{k}^{\varepsilon,t,x,*}=-(M_{k,2}^\varepsilon)^{-1}M_{k,1}
(X_{k}^{\varepsilon,t,x,*}, u^{\varepsilon,t,x,*}|_{\mathbb{T}_{k+1}}),~~
k\in \mathbb{T}_t,
\end{eqnarray*}
and
\begin{eqnarray}\label{open-loop-equilibrium-epsilon}
\left\{
\begin{array}{l}
X^{\varepsilon,t,x,*}_{k+1} =\big{[}\mathcal{A}_{k,k}
X^{\varepsilon,t,x,*}_{k}+\mathcal{B}_{k,k}u^{\varepsilon,t,x,*}_k
+f_{k,k}\big{]}\\[1mm]
\hphantom{X^{\varepsilon,t,x,*}_{k+1}}+\big{[}\mathcal{C}_{k,k}X^{\varepsilon,t,x,*}_{k}
+\mathcal{D}_{k,k}u^{\varepsilon,t,x,*}_k+d_{k,k}\big{]}w_k,\\[1mm]
X^{\varepsilon,t,x,*}_{t} = x,~~ k \in  \mathbb{T}_t.
\end{array}
\right.
\end{eqnarray}
By (\ref{open-loop-equilibrium-k}), $u^{\varepsilon, t, x, *}$ can
be expressed as
\begin{eqnarray}\label{open-loop-equilibrium-k-3}
u^{\varepsilon, t,x,*}_k=-(\mathcal{W}_{k}^{\varepsilon})^\dagger \mathcal{H}^\varepsilon_{k} X_{k}^{\varepsilon, t,x,*}-(\mathcal{W}_{k}^{\varepsilon})^\dagger \beta^\varepsilon_{k},~~~~k\in \mathbb{T}_t.
\end{eqnarray}
Here, $\mathcal{W}_{k}^{\varepsilon}, \mathcal{H}^\varepsilon_{k}$
and $\beta^\varepsilon_{k}$ are obtained by replacing
$\mathcal{R}_{k,k}$ with $\mathcal{R}_{k,k}+\varepsilon I$ in
(\ref{feedback-gain-k}).
%
%
%
%
%
%
%
%
%

\begin{theorem}\label{Theorem-epsilon}
Let
\begin{eqnarray}\label{p-3}
\left\{
\begin{array}{l}
\left\{
\begin{array}{l}
P_{k,\ell}=Q_{k,\ell}+A^T_{k,\ell}P_{k,\ell+1}A_{k,\ell}+C^T_{k,\ell}P_{k,\ell+1}C_{k,\ell},\\[1mm]
{\mathcal{P}}_{k,\ell}=\mathcal{Q}_{k,\ell}+\mathcal{A}^T_{k,\ell}\mathcal{P}_{k,\ell+1}{\mathcal{A}}_{k,\ell}+\mathcal{C}^T_{k,\ell}P_{k,\ell+1}{\mathcal{C}}_{k,\ell},\\[1mm]
P_{k,N}=G_k,~~\mathcal{{P}}_{k,N}=\mathcal{{G}}_k,~~\ell\in \mathbb{T}_k,
\end{array}
\right. \\[1mm]
\mathcal{R}_{k,k}+\mathcal{B}_{k,k}^T\mathcal{P}_{k,k+1}\mathcal{B}_{k,k}
+\mathcal{D}_{k,k}^TP_{k,k+1}\mathcal{D}_{k,k}\geq 0,\\[1mm]
k\in \mathbb{T}
\end{array}
\right.
\end{eqnarray}
%
%
%
%
%
be solvable. Then the following statements hold.

(i) For any $t\in \mathbb{T}$ and any $x\in L^2_{\mathcal{F}}(t;
\mathbb{R}^n)$,  Problem (LQ)$_\varepsilon$ for the initial pair
$(t,x)$ admits an open-loop equilibrium pair
$(X^{\varepsilon,t,x,*}, u^{\varepsilon,t,x,*})$, which is given
in (\ref{open-loop-equilibrium-epsilon}) and
(\ref{open-loop-equilibrium-k-3}). 

(ii) If the sequence $\{u^{\varepsilon,t,x,*}, \varepsilon>0\}$ is
bounded in  $L^2_{\mathcal{F}}(\mathbb{T}_t; \mathbb{R}^m)$, then
Problem (LQ) for the initial pair $(t,x)$ admits an open-loop
equilibrium pair.

\end{theorem}

\emph{Proof}. (i) follows directly from the comments above.

(ii). As $\{u^{\varepsilon,t,x,*}, \varepsilon>0\}$ is bounded in $L^2_{\mathcal{F}}(\mathbb{T}_t; \mathbb{R}^m)$, there exists a weakly convergent subsequence $\{u^{\varepsilon_j,t,x,*}, j\in \{0,1,2,...\}\}$ of
$\{u^{\varepsilon,t,x,*}, \varepsilon>0\}$ with its weak limit ${\overline{v}}^{t,x,*}$. We can further select a subsequence of $\{u^{\varepsilon_j,t,x,*}, j\in \{0,1,2,...\}\}$ such that the subsequence converges strongly to ${\overline{v}}^{t,x,*}$. Without loss of generality, we assume that $\{u^{\varepsilon_j,t,x,*}, j\in \{0,1,2,...\}\}$ converges to ${\overline{v}}^{t,x,*}$ strongly. Denote $\overline{X}^{t,x,*}$ as the solution to the following equation
\begin{eqnarray}\label{open-loop-equilibrium-v}
\left\{
\begin{array}{l}
\overline{X}^{t,x,*}_{k+1} =\big{[}\mathcal{A}_{k,k}\overline{X}^{t,x,*}_{k}+\mathcal{B}_{k,k}\overline{v}^{t,x,*}_k+f_{k,k}\big{]}+\big{[}\mathcal{C}_{k,k}\overline{X}^{t,x,*}_{k} +\mathcal{D}_{k,k}\overline{v}^{t,x,*}_k+d_{k,k}\big{]}w_k,\\[1mm]
\overline{X}^{t,x,*}_{t} = x,~~ k \in  \mathbb{T}_t.
\end{array}
\right.
\end{eqnarray}
%
Clearly,
\begin{eqnarray}\label{open-loop-equilibrium-v-2}
\overline{X}^{t,x,*}_k=\widehat{X}^{t,x,*}_k+\widetilde{X}^{t,0,*}_k,~k \in  \mathbb{T}_t,
\end{eqnarray}
where $\widehat{X}^{t,x,*}_k$ and $\widetilde{X}^{t,0,*}_k$ are computed via
\begin{eqnarray*}
\left\{
\begin{array}{l}
\widehat{X}^{t,x,*}_{k+1} =\big{[}\mathcal{A}_{k,k}\widehat{X}^{t,x,*}_{k}+f_{k,k}\big{]}+\big{[}\mathcal{C}_{k,k}\widehat{X}^{t,x,*}_{k} +d_{k,k}\big{]}w_k,\\[1mm]
\widehat{X}^{t,x,*}_{t} = x,~~ k \in  \mathbb{T}_t,
\end{array}
\right.
\end{eqnarray*}
and
\begin{eqnarray*}
\left\{\begin{array}{l}
\widetilde{{X}}^{t,0,*}_{k+1} =\big{[}\mathcal{A}_{k,k}\widetilde{X}^{t,0,*}_{k}+\mathcal{B}_{k,k}\overline{v}^{t,x,*}_k\big{]} +\big{[}\mathcal{C}_{k,k}\widetilde{X}^{t,0,*}_{k} +\mathcal{D}_{k,k}\overline{v}^{t,x,*}_k\big{]}w_k,\\[1mm]
\widetilde{X}^{t,0,*}_{t} = 0,~~ k \in  \mathbb{T}_t.
\end{array}
\right.
\end{eqnarray*}
From (\ref{open-loop-equilibrium-v-2}), (\ref{open-loop-equilibrium-v}) actually introduces an affine
operator, which is defined from
$L^2_{\mathcal{F}}(\mathbb{T}_t; \mathbb{R}^m)$ to
$L^2_{\mathcal{F}}(\widetilde{\mathbb{T}}_t; \mathbb{R}^n)$, i.e.,
$\overline{X}^{t,x,*}=\widehat{X}^{t,x,*}+\Theta(\overline{v}^{t,x,*})$ with $\widetilde{X}^{t,0,*}=\Theta(\overline{v}^{t,x,*})$. It can be shown
that the operator $\Theta$ is linear and bounded. As
$u^{\varepsilon_j,t,x,*}$ converges strongly to
${\overline{v}}^{t,x,*}$, ${X}^{\varepsilon,t,x,*}=\widehat{X}^{t,x,*}+\Theta(u^{\varepsilon,t,x,*})$ will converge
strongly to $\overline{X}^{t,x,*}=\widehat{X}^{t,x,*}+\Theta(\overline{v}^{t,x,*})$.
Furthermore, from the definition of open-loop equilibrium control, it
follows that for any $k\in \mathbb{T}_t$ and any $u_k\in L^2_\mathcal{F}(k;
\mathbb{R}^m)$
\begin{eqnarray}\label{open-loop-equilibrium-k-6}
&&\hspace{-3em}J(k,X^{\varepsilon_j,t,x,*}_k; u^{\varepsilon_j,t,x,*}|_{\mathbb{T}_k})
+\varepsilon_j \mathbb{E}|u^{\varepsilon_j,t,x,*}_k|^2\nonumber \\[1mm]
&&\hspace{-3em}\leq J(k,X^{\varepsilon_j,t,x,*}_k;(u_k, u^{\varepsilon_j,t,x,*}
|_{\mathbb{T}_{k+1}}))+\varepsilon_j \mathbb{E}|u_k|^2.
\end{eqnarray}
Letting $j\rightarrow \infty$ in (\ref{open-loop-equilibrium-k-6}), we have for any $\forall k\in \mathbb{T}_t$ and any $u_k\in L^2_\mathcal{F}(k; \mathbb{R}^m)$
\begin{eqnarray}\label{open-loop-equilibrium-k-7}
J(k,\overline{X}^{t,x,*}_k; \overline{v}^{t,x,*}|_{\mathbb{T}_k})
\leq J(k,\overline{X}^{t,x,*}_k;(u_k, \overline{v}^{t,x,*}|_{\mathbb{T}_{k+1}})).
\end{eqnarray}
By (\ref{open-loop-equilibrium-v}), we can assert that
$(\overline{X}^{t,x,*}, \overline{v}^{t,x,*})$ is an open-loop
equilibrium pair of Problem (LQ) for the initial pair $(t,x)$.\endpf

Under the point-wise convexity condition (\ref{p-3}), above theorem tells us that Problem (LQ) is ``almostly" solvable: for any arbitrarily small $\varepsilon$, the perturbation version Problem (LQ)$_\varepsilon$ of Problem (LQ) is solvable for any initial pair. Naturally, we may ask: when will Problem (LQ) be solvable for all the initial pairs? The following lemma presents a sufficient condition to the existence of the open-loop
equilibrium control of Problem (LQ).

\begin{lemma}\label{Lemma-sufficiency}
For $\mathcal{W}_k, \mathcal{H}_k, \beta_k, k\in \mathbb{T}_t$ (defined in (\ref{feedback-gain-k})), if
\begin{eqnarray}\label{W-H-beta}
\mathcal{W}_k\mathcal{W}_k^\dagger\mathcal{H}_k-\mathcal{H}_k=0,~~
\mathcal{W}_k\mathcal{W}_k^\dagger\beta_k-\beta_k=0,~~k\in\mathbb{T}_t
\end{eqnarray}
are satisfied and (\ref{P}) is solvable, then Problem (LQ) for the initial pair
$(t,x)$ admits an open-loop equilibrium pair.
\end{lemma}

\emph{Proof}. Introduce a dynamics
\begin{eqnarray}\label{open-loop-equilibrium-state-suffi}
\left\{
\begin{array}{l}
\widetilde{X}^{t,x,*}_{k+1} = \big{[}\big{(}\mathcal{A}_{k,k}
-\mathcal{B}_{k,k}\mathcal{W}^\dagger_k\mathcal{H}_k\big{)}
\widetilde{X}^{t,x,*}_{k}+f_{k,k}-\mathcal{B}_{k,k}\mathcal{W}^\dagger_k\beta_k\big{]}\\[1mm]
\hphantom{X^{t,x,*}_{k+1}=}+\big{[}\big{(}\mathcal{A}_{k,k}
-\mathcal{B}_{k,k}\mathcal{W}^\dagger_k\mathcal{H}_k\big{)}
\widetilde{X}^{t,x,*}_{k}+d_{k,k}-\mathcal{D}_{k,k}\mathcal{W}^\dagger_k\beta_k\big{]}w_k,\\[1mm]
\widetilde{X}^{t,x,*}_{t} = x,~~ k \in  \mathbb{T}_t,
\end{array}
\right.
\end{eqnarray}
and a control
\begin{eqnarray}\label{open-loop-equilibrium-k-suffi}
\widetilde{u}_{k}^{t,x,*}=-\mathcal{W}_{k}^\dagger \mathcal{H}_{k}
\widetilde{X}_{k}^{t,x,*}-\mathcal{W}_{k}^\dagger \beta_{k},~~k\in
\mathbb{T}_t.
\end{eqnarray}
Then, by reversing the first part of the proof of Theorem
\ref{Theorem-Necessary}, we can show that for any $k\in
\mathbb{T}_t$, the following FBS$\Delta$E admits a solution
\begin{eqnarray*}
\left\{\begin{array}{l}
\widetilde{X}^{k,t,x}_{\ell+1}=\big{(}A_{k,\ell}\widetilde{X}^{k,t,x}_\ell
+\bar{A}_{k,\ell}\mathbb{E}_k\widetilde{X}^{k,t,x}_\ell+B_{k,\ell}\widetilde{u}^{t,x,*}_\ell+\bar{B}_{k,\ell}
\mathbb{E}_k\widetilde{u}^{t,x,*}_\ell+f_{k,\ell}\big{)}\\[1mm]
\hphantom{\widetilde{X}^{k,t,x}_{\ell+1}=}+\big{(}C_{k,\ell}\widetilde{X}^{k,t,x}_\ell
+\bar{C}_{k,\ell}\mathbb{E}_k\widetilde{X}^{k,t,x}_\ell+D_{k,\ell}\widetilde{u}^{t,x,*}_\ell
+\bar{D}_{k,\ell}\mathbb{E}_k\widetilde{u}^{t,x,*}_\ell+d_{k,\ell}\big{)}w_\ell, \\[1mm]
\widetilde{Z}_\ell^{k,t,x}=A_{k,\ell}^T\mathbb{E}_\ell \widetilde{Z}_{\ell+1}^{k,t,x}
+ \bar{A}_{k,\ell}^T\mathbb{E}_k\widetilde{Z}_{\ell+1}^{k,t,x}+C_{k,\ell}^T\mathbb{E}_\ell(\widetilde{Z}_{\ell+1}^{k,t,x}w_\ell)
+\bar{C}_{k,\ell}^T\mathbb{E}_k(\widetilde{Z}_{\ell+1}^{k,t,x}w_\ell)\\[1mm]
\hphantom{\widetilde{Z}_\ell^{k,t,x}=}+Q_{k,\ell}\widetilde{X}_\ell^{k,t,x}
+\bar{Q}_{k,\ell}\mathbb{E}_k\widetilde{X}_\ell^{k,t,x}+q_{k,\ell},\\[1mm]
\widetilde{X}^{k,t,x}_k=\widetilde{X}^{t,x,*}_k,~~\\[1mm]
\widetilde{Z}_N^{k,t,x}
=G_k\widetilde{X}^{k,t,x}_N+\bar{G}_k\mathbb{E}_k\widetilde{X}^{k,t,x}_N+g_k,~~
\ell\in \mathbb{T}_k
\end{array}
\right.
\end{eqnarray*}
with properties
\begin{eqnarray*}
&&\widetilde{Z}^{k,t,x}_\ell=P_{k,\ell}\widetilde{X}^{k,t,x}_\ell+\bar{P}_{k,\ell}
\mathbb{E}_k\widetilde{X}^{k,t,x}_\ell+T_{k,\ell}\widetilde{X}^{t,x,*}_\ell
+\bar{T}_{k,\ell}\mathbb{E}_k\widetilde{X}^{t,x,*}_\ell+\pi_{k,\ell},~~
\ell\in \mathbb{T}_{k},
\end{eqnarray*}
and
\begin{eqnarray*}
0=\mathcal{R}_{k,k}\widetilde{u}_k^{t,x,*}+\mathcal{B}_{k,k}^T\mathbb{E}_k
\widetilde{Z}^{k,t,x}_{k+1}
+\mathcal{D}_{k,k}^T\mathbb{E}_k(\widetilde{Z}_{k+1}^{k,t,x}w_k)+\rho_{k,k}.
\end{eqnarray*}
Furthermore, by (\ref{convex-4}) and (\ref{convex-5}) we have
(\ref{convex}). From Theorem \ref{Theorem-Equivalentce-open-loop},
$(\widetilde{X}^{t,x,*},
\widetilde{u}^{t,x,*})$ is an open-loop equilibrium pair of Problem (LQ) for the initial pair $(t,x)$. This completes the
proof. \endpf

In fact, we have the following necessary and sufficient condition to the existence of open-loop equilibrium pair \emph{for any initial pair}.

\begin{theorem}\label{Theorem-Necessary-sufficient}
The following statements are equivalent.

(i) For any $t\in \mathbb{T}$ and any $x\in L^2_\mathcal{F}(t;
\mathbb{R}^n)$, Problem (LQ) admits an
open-loop equilibrium pair for the initial pair $(t,x)$.

(ii) (\ref{p-3}), the set of GDREs
\begin{eqnarray}\label{T-3}
\left\{
\begin{array}{l}
\left\{
\begin{array}{l}
T_{k,\ell}=A^T_{k,\ell}T_{k,\ell+1}\mathcal{A}_{\ell,\ell}
+C^T_{k,\ell}T_{k,\ell+1}\mathcal{C}_{\ell,\ell}\\[1mm]
\hphantom{T_{k,\ell}=}-\Big{(}A^T_{k,\ell}P_{k,\ell+1}B_{k,\ell}
+A^T_{k,\ell}T_{k,\ell+1}\mathcal{B}_{\ell,\ell}+C^T_{k,\ell}P_{k,\ell+1}D_{k,\ell}+C^T_{k,\ell}T_{k,\ell+1}
\mathcal{D}_{\ell,\ell}\Big{)}\mathcal{W}_{\ell}^\dagger \mathcal{H}_{\ell},\\[1mm]
{\mathcal{T}}_{k,\ell}=\mathcal{A}^T_{k,\ell}{\mathcal{T}}_{k,\ell+1}
\mathcal{A}_{\ell,\ell}+{\mathcal{C}}^T_{k,\ell}T_{k,\ell+1}
\mathcal{C}_{\ell,\ell}\\[1mm]
\hphantom{\bar{T}_{k,\ell}=}-\Big{(}\mathcal{A}^T_{k,\ell}
\mathcal{P}_{k,\ell+1}{\mathcal{B}}_{k,\ell}
+\mathcal{A}^T_{k,\ell}{\mathcal{T}}_{k,\ell+1}\mathcal{B}_{\ell,\ell}+\mathcal{C}^T_{k,\ell}P_{k,\ell+1}{\mathcal{D}}_{k,\ell}+{\mathcal{C}}^T_{k,\ell}
T_{k,\ell+1}\mathcal{D}_{\ell,\ell}\Big{)}\mathcal{W}_{\ell}^\dagger \mathcal{H}_{\ell},\\[1mm]
T_{k,N}=0,~~{\mathcal{T}}_{k,N}=0, \\[1mm]
\ell\in \mathbb{T}_{k},
\end{array}
\right.\\[1mm]
\mathcal{W}_k\mathcal{W}_k^\dagger\mathcal{H}_k-\mathcal{H}_k=0,\\[1mm]
k\in \mathbb{T},
\end{array}
\right.
\end{eqnarray}
and the set of LDEs
\begin{eqnarray}\label{pi-3}
\left\{
\begin{array}{l}
\left\{
\begin{array}{l}
\pi_{k,\ell}=\mathcal{A}^T_{k,\ell}\mathcal{P}_{k,\ell+1}\big{(}f_{k,\ell}
-\mathcal{B}_{k,\ell}\mathcal{W}_{\ell}^\dagger \beta_{\ell}\big{)}+\mathcal{A}^T_{k,\ell}\mathcal{T}_{k,\ell+1}\big{(}f_{\ell,\ell}
-\mathcal{B}_{\ell,\ell} \mathcal{W}_{\ell}^\dagger \beta_{\ell}\big{)}
\\[1mm]
\hphantom{\pi_{k,\ell}=}+\mathcal{C}^T_{k,\ell}P_{k,\ell+1}\big{(}d_{k,\ell}
-\mathcal{D}_{k,\ell}\mathcal{W}_{\ell}^\dagger \beta_{\ell}\big{)}+\mathcal{C}^T_{k,\ell}T_{k,\ell+1}\big{(}d_{\ell,\ell}-\mathcal{D}_{\ell,\ell}
\mathcal{W}_{\ell}^\dagger \beta_{\ell}\big{)}\\[1mm]
\hphantom{\pi_{k,\ell}=}+\mathcal{A}^T_{k,\ell}\pi_{k,\ell+1}+q_{k,\ell},\\[1mm]
\pi_{k,N}=g_k,
\end{array}
\right.\\[1mm]
\mathcal{W}_k\mathcal{W}_k^\dagger\beta_k-\beta_k=0,\\[1mm]
k\in \mathbb{T}_t
\end{array}
\right.
\end{eqnarray}
are solvable in the sense
\begin{eqnarray*}
\left\{
\begin{array}{l}
\mathcal{R}_{k,k}+\mathcal{B}_{k,k}^T\mathcal{P}_{k,k+1}\mathcal{B}_{k,k}
+\mathcal{D}_{k,k}^TP_{k,k+1}\mathcal{D}_{k,k}\geq 0,\\[1mm]
\mathcal{W}_k\mathcal{W}_k^\dagger\mathcal{H}_k-\mathcal{H}_k=0,\\[1mm]
\mathcal{W}_k\mathcal{W}_k^\dagger\beta_k-\beta_k=0,\\[1mm]
k\in \mathbb{T}.
\end{array}
\right.
\end{eqnarray*}
Here,
\begin{eqnarray*}\label{feedback-gain-k-3}
\left\{
\begin{array}{l}
\mathcal{W}_{k}=\mathcal{R}_{k,k}+\mathcal{B}_{k,k}^T\big{(}\mathcal{P}_{k,k+1}+\mathcal{T}_{k,k+1}\big{)}
\mathcal{B}_{k,k}+\mathcal{D}_{k,k}^T\big{(}{P}_{k,k+1}+T_{k,k+1}\big{)}\mathcal{D}_{k,k},\\[1mm]
\mathcal{H}_{k}=\mathcal{B}_{k,k}^T\big{(}\mathcal{P}_{k,k+1}+\mathcal{T}_{k,k+1}\big{)}
\mathcal{A}_{k,k}+\mathcal{D}_{k,k}^T\big{(}{P}_{k,k+1}+T_{k,k+1}\big{)}\mathcal{C}_{k,k},\\[1mm]
\beta_{k}=\mathcal{B}_{k,k}^T\big{[}\big{(}\mathcal{{P}}_{k,k+1}+\mathcal{T}_{k,k+1}\big{)}
f_{k,k}+\pi_{k,k+1}\big{]}+\mathcal{D}_{k,k}^T\big{(}{P}_{k,k+1}+T_{k,k+1}\big{)}d_{k,k}
+\rho_{k,k},\\[1mm] k\in \mathbb{T}.
\end{array}
\right.
\end{eqnarray*}

Under any of above conditions, an open-loop equilibrium control for the initial pair
$(t,x)$ is given in (\ref{open-loop-equilibrium-k}).

\end{theorem}

\emph{Proof}. The sufficiency follows from Lemma
\ref{Lemma-sufficiency}. As for the necessity, by Theorem
\ref{Theorem-Necessary} we need only to prove
\begin{eqnarray}\label{W-H-3}
\mathcal{W}_k\mathcal{W}_k^\dagger\mathcal{H}_k-\mathcal{H}_k=0,~~
\mathcal{W}_k\mathcal{W}_k^\dagger\beta_k-\beta_k=0,~~k\in\mathbb{T}.
\end{eqnarray}
Consider Problem (LQ) for the initial pair $(N-1,x)$ with $x\in
L^2_\mathcal{F}(N-1; \mathbb{R}^n)$. By the proof of Theorem
\ref{Theorem-Necessary}, we have
\begin{eqnarray}\label{N-1}
0=\mathcal{W}_{N-1}u_{N-1}^{N-1,x,*}+\mathcal{H}_{N-1}
X_{N-1}^{N-1,x,*}+\beta_{N-1}.
\end{eqnarray}
Noting $X_{N-1}^{N-1,x,*}=x$ and taking  $x=0$ in (\ref{N-1}), we
have
\begin{eqnarray}\label{jif-a1}
0=\mathcal{W}_{N-1}u_{N-1}^{N-1, 0,*}+ \beta_{N-1},
\end{eqnarray}
which together with Lemma \ref{Lemma-matrix-equation} leads to
$\mathcal{W}_{N-1}\mathcal{W}_{N-1}^\dagger\beta_{N-1}-\beta_{N-1}
=0$.
Furthermore, by subtracting (\ref{jif-a1}) from (\ref{N-1}) we have
\begin{eqnarray*}
0=\mathcal{W}_{N-1}\big{(}u_{N-1}^{N-1, x,*}-u_{N-1}^{N-1, 0,*}\big{)}+\mathcal{H}_{N-1} x.
\end{eqnarray*}
Let $e_i$ be a $\mathbb{R}^n$-valued vector with its $i$-th entry
being 1 and other entries 0. Then, we have
\begin{eqnarray*}
&&\hspace{-2em}0=\mathcal{W}_{N-1}\big{(}u_{N-1}^{N-1, e_1,*}-u_{N-1}^{N-1,
0,*},...,u_{N-1}^{N-1, e_n,*}-u_{N-1}^{N-1,
0,*}\big{)}+\mathcal{H}_{N-1}\big{(}e_1,...,e_n \big{)}.
\end{eqnarray*}
Noting that $\big{(}e_1,...,e_n \big{)}$ is the identity matrix and
by Lemma \ref{Lemma-matrix-equation}, we have
$\mathcal{W}_{N-1}\mathcal{W}_{N-1}^\dagger\mathcal{H}_{N-1}-\mathcal{H}_{N-1}=0$.

Considering Problem (LQ) for the initial pair $(N-2,x)$ with $x\in
L^2_\mathcal{F}(N-2; \mathbb{R}^n)$, we can similarly prove
\begin{eqnarray*}
&&\mathcal{W}_{N-2}\mathcal{W}_{N-2}^\dagger\mathcal{H}_{N-2}-\mathcal{H}_{N-2}=0,~~\mathcal{W}_{N-2}\mathcal{W}_{N-2}^\dagger\beta_{N-2}-\beta_{N-2}=0.
\end{eqnarray*}
Continuing above procedure backwardly, we then achieve the conclusion.\endpf 

Note that $P_{k,\ell}, \mathcal{P}_{k,\ell}, k\in \mathbb{T},
\ell\in \mathbb{T}_{k+1}$, are symmetric. If $Q_{k,\ell},
\bar{Q}_{k,\ell}, R_{k,\ell}, \bar{R}_{k,\ell}$ are selected such
that
$
Q_{k,\ell},~ Q_{k,\ell}+\bar{Q}_{k,\ell}, ~R_{k,\ell},~
R_{k,\ell}+\bar{R}_{k,\ell}\geq 0, ~k\in \mathbb{T}, ~ \ell\in
\mathbb{T}_k,
$
then (\ref{p-3}) is solvable. Furthermore, as indicated in
(\ref{Z}), $\Theta_{k,\ell}=\{P_{k,\ell}, \mathcal{P}_{k,\ell},
T_{k,\ell}, \mathcal{T}_{k,\ell},\pi_{k,\ell}\}$ is used to express
$Z^{k,t,x}_{\ell}$. $\{P_{k,\ell}, \mathcal{P}_{k,\ell}\}$ is then
called the symmetric part of $\Theta_{k,\ell}$, and $\{T_{k,\ell},
\mathcal{T}_{k,\ell}\}$, $\pi_{k,\ell}$ are viewed as the
nonsymmetric part and nonhomogeneous part, respectively.

To end this section, we give some comments on the open-loop equilibrium control $u^{t,x,*}$ and its closed-loop expression  (\ref{open-loop-equilibrium-k}). Generally speaking, for deterministic problems, an open-loop control is a functional of initial state and the time variable, and a closed-loop control is a functional of the observed state information. As all the states are essentially functionals of initial state and the time variable, a closed-loop control is an open-loop control indeed. Concerned with the stochastic case, a control problem is formulated within a random background, which is characterized via a filtration.
%
%
It is better to select the open-loop control to be adapted to the background filtration. For example, an open-loop control in this paper is selected to be adapted to $\{\mathcal{F}_k\}$. The closed-loop control is similarly defined as that for the deterministic case, and also a closed-loop control is an open-loop control. Therefore, though (\ref{open-loop-equilibrium-k}) is of closed-loop form, it is indeed an open-loop control.

\section{Some Special Cases}\label{section-3}

\subsection{The state matrices and weighting matrices are independent of the initial time}

In this case, the system dynamics and the cost functional are,
respectively, given by
\begin{eqnarray}\label{system-1-mf}
\left\{\begin{array}{l}
X^t_{k+1}=\big{(}A_{k}X^t_k+\bar{A}_{k}\mathbb{E}_tX^t_k+B_{k}u_k+\bar{B}_{k}
\mathbb{E}_tu_k+f_{k}\big{)}\\[1mm]
\hphantom{X^t_{k+1}=}+\big{(}C_{k}X^t_k+\bar{C}_{k}\mathbb{E}_tX^t_k+D_{k}u_k
+\bar{D}_{k}\mathbb{E}_tu_k+d_{k}\big{)}w_k, \\[1mm]
X^t_t=x,~~k\in \mathbb{T}_t,~~t\in \mathbb{T},
\end{array}
\right.
\end{eqnarray}
and
\begin{eqnarray}\label{cost-1-mf}
&&\hspace{-2em}J(t,x;u)
 =\sum_{k=t}^{N-1}\dbE_t\Big{[}(X_k^t)^TQ_{k}X^t_k +(\mathbb{E}_tX_k^t)^T
 \bar{Q}_{k}\mathbb{E}_tX^t_k+ u_k^TR_{k}u_k+(\mathbb{E}_tu_k)^T \bar{R}_{k}
 \mathbb{E}_tu_k  +2q_{k}^T X^t_k+ 2\rho_{k}^T u_k\Big{]}\nonumber \\[1mm]
&&\hspace{-2em}
 \hphantom{J(t,x;u)
 =}
+\mathbb{E}_t\big{[}(X_N^t)^TGX^t_N\big{]}
+(\mathbb{E}_tX^t_N)^T \bar{G}\mathbb{E}_tX^t_{N}+2\mathbb{E}_tg^T X^t_N.
\end{eqnarray}
%
Problem (LQ) corresponding to (\ref{system-1-mf}) and
(\ref{cost-1-mf}) will be denoted as Problem (LQ)$_{s1}$. Now,
(\ref{p-3}), (\ref{T-3}) and (\ref{pi-3}) become
\begin{eqnarray}\label{P-1-mf}
\left\{
\begin{array}{l}
P_{k}=Q_{k}+A^T_{k}P_{k+1}A_{k}+C^T_{k}P_{k+1}C_{k},\\[1mm]
{\mathcal{P}}_{k}=\mathcal{Q}_{k}+\mathcal{A}^T_{k}
\mathcal{P}_{k+1}{\mathcal{A}}_{k}+\mathcal{C}^T_{k}P_{k+1}{\mathcal{C}}_{k},\\[1mm]
P_{N}=G,~~{\mathcal{P}}_{N}={\mathcal{G}},
\\[1mm]
\mathcal{R}_{k}+\mathcal{B}_{k}^T\mathcal{P}_{k+1}\mathcal{B}_{k}
+\mathcal{D}_{k}^TP_{k+1}\mathcal{D}_{k}\geq 0,\\[1mm]
k\in \mathbb{T},
\end{array}
\right.
\end{eqnarray}
\begin{eqnarray}\label{T-1-mf}
\left\{
\begin{array}{l}
T_{k}=A^T_{k}T_{k+1}\mathcal{A}_{k}+C^T_{k}T_{k+1}\mathcal{C}_{k}\\[1mm]
\hphantom{T_{k}=}-\big{(}A^T_{k}P_{k+1}B_{k}+A^T_{k}T_{k+1}\mathcal{B}_{k}
+C^T_{k}P_{k+1}D_{k}+C^T_{k}T_{k+1}\mathcal{D}_{k}\big{)}
\mathcal{W}_{k}^\dagger \mathcal{H}_{k},\\[1mm]
{\mathcal{T}}_{k}=\mathcal{A}^T_{k}{\mathcal{T}}_{k+1}
\mathcal{A}_{k}+{\mathcal{C}}^T_{k}T_{k+1}\mathcal{C}_{k}\\[1mm]
\hphantom{{\mathcal{T}}_{k}=}-\Big{(}
\mathcal{A}^T_{k}{(}\mathcal{P}_{k+1}+{\mathcal{T}}_{k+1}{)}
\mathcal{B}_{k}+\mathcal{C}^T_{k}{(}P_{k+1}+T_{k+1}{)}
\mathcal{D}_{k}\Big{)}\mathcal{W}_{k}^\dagger \mathcal{H}_{k},\\[1mm]
T_{N}=0,~~{\mathcal{T}}_{N}=0, \\[1mm]
\mathcal{W}_k\mathcal{W}_k^\dagger\mathcal{H}_k-\mathcal{H}_k=0,\\[1mm]
k\in \mathbb{T},
\end{array}
\right.
\end{eqnarray}
and
\begin{eqnarray}\label{pi-1-mf}
\left\{
\begin{array}{l}
\pi_{k}=\mathcal{A}^T_{k}\mathcal{P}_{k+1}\big{(}f_{k}-\mathcal{B}_{k}
\mathcal{W}_{k}^\dagger \beta_{k}\big{)}+\mathcal{A}^T_{k}
\mathcal{T}_{k+1}\big{(}f_{k}-\mathcal{B}_{k} \mathcal{W}_{k}^\dagger
\beta_{k}\big{)}+\mathcal{A}^T_{k}\pi_{k+1}\\[1mm]
\hphantom{\pi_{k}=}+\mathcal{C}^T_{k}P_{k+1}\big{(}d_{k}
-\mathcal{D}_{k}\mathcal{W}_{k}^\dagger \beta_{k}\big{)}+\mathcal{C}^T_{k}T_{k+1}\big{(}d_{k}-\mathcal{D}_{k}
\mathcal{W}_{k}^\dagger \beta_{k}\big{)}+q_{k},\\[1mm]
\pi_{N}=g,\\[1mm]
\mathcal{W}_k\mathcal{W}_k^\dagger\beta_k-\beta_k=0,\\[1mm]
k\in \mathbb{T},
\end{array}
\right.
\end{eqnarray}
where
\begin{eqnarray*}\label{feedback-gain-k-mf}
\left\{
\begin{array}{l}
\mathcal{W}_{k}=\mathcal{R}_{k}+\mathcal{B}_{k}^T\big{(}\mathcal{P}_{k+1}
+\mathcal{T}_{k+1}\big{)}\mathcal{B}_{k}+\mathcal{D}_{k}^T\big{(}{P}_{k+1}+
T_{k+1}\big{)}\mathcal{D}_{k},\\[1mm]
\mathcal{H}_{k}=\mathcal{B}_{k}^T\big{(}\mathcal{P}_{k+1}+\mathcal{T}_{k+1}
\big{)}\mathcal{A}_{k}+\mathcal{D}_{k}^T\big{(}{P}_{k+1}+T_{k+1}\big{)}
\mathcal{C}_{k},\\[1mm]
\beta_{k}=\mathcal{B}_{k}^T\big{[}\big{(}\mathcal{{P}}_{k+1}
+\mathcal{T}_{k+1}\big{)}f_{k}+\pi_{k+1}\big{]}+\mathcal{D}_{k}^T
\big{(}{P}_{k+1}+T_{k+1}\big{)}d_{k}+\rho_{k}.\\[1mm]
k\in \mathbb{T}.
\end{array}
\right.
\end{eqnarray*}
%
By Theorem
\ref{Theorem-Necessary-sufficient}, we have the following result.

\begin{corollary}\label{corollary-Necessary-sufficient}
For any $t\in \mathbb{T}$ and any $x\in L^2_\mathcal{F}(t;
\mathbb{R}^n)$, Problem (LQ)$_{s1}$ for the initial pair $(t,x)$
admits an open-loop equilibrium pair if and only if (\ref{P-1-mf}), (\ref{T-1-mf}) and (\ref{pi-1-mf}) are solvable.
\end{corollary}

\subsection{The case without mean-field terms}

Consider the following system dynamics and cost functional
\begin{eqnarray}\label{system-1-in}
\left\{\begin{array}{l}
X^t_{k+1}=\big{(}A_{t,k}X^t_k+B_{t,k}u_k+f_{t,k}\big{)}+\big{(}C_{t,k}X^t_k+D_{t,k}u_k+d_{t,k}\big{)}w_k, \\[1mm]
X^t_t=x,~~k\in \mathbb{T}_t,~~t\in \mathbb{T},
\end{array}
\right.
\end{eqnarray}
and
\begin{eqnarray}\label{cost-1-in}
&&\hspace{-2em}J(t,x;u)
 =\sum_{k=t}^{N-1}\dbE_t\big{[}(X_k^t)^TQ_{t,k}X^t_k + u_k^TR_{t,k}u_k+2q_{t,k}^T X^t_k+ 2\rho_{t,k}^T u_k\big{]}\nonumber \\[1mm]
&&\hspace{-2em}\hphantom{J(t,x;u)
 =}
+\mathbb{E}_t\big{[}(X_N^t)^TG_tX^t_N\big{]}+2\mathbb{E}_tg_t^T X^t_N. 
\end{eqnarray}
Problem (LQ) corresponding to (\ref{system-1-in}) and
(\ref{cost-1-in}) will be denoted as Problem (LQ)$_{s2}$. In this
case, we have
\begin{eqnarray}\label{P-2-s2}
\left\{
\begin{array}{l}
\left\{
\begin{array}{l}
P_{k,\ell}=Q_{k,\ell}+A^T_{k,\ell}P_{k,\ell+1}A_{k,\ell}+C^T_{k,\ell}P_{k,\ell+1}C_{k,\ell},\\[1mm]
P_{k,N}=G_k,
~~\ell\in \mathbb{T}_k,
\end{array}
\right. \\[1mm]
{R}_{k,k}+{B}_{k,k}^T{P}_{k,k+1}{B}_{k,k}+{D}_{k,k}^TP_{k,k+1}{D}_{k,k}\geq 0,\\[1mm]
k\in \mathbb{T},
\end{array}
\right.
\end{eqnarray}
%
%
\begin{eqnarray}\label{T-2-s2}
\left\{
\begin{array}{l}
\left\{
\begin{array}{l}
T_{k,\ell}=A^T_{k,\ell}T_{k,\ell+1}{A}_{\ell,\ell}+C^T_{k,\ell}T_{k,\ell+1}{C}_{\ell,\ell}\\[1mm]
\hphantom{T_{k,\ell}=}-\big{(}A^T_{k,\ell}(P_{k,\ell+1}+T_{k,\ell+1}){B}_{\ell,\ell}
+C^T_{k,\ell}(P_{k,\ell+1}+T_{k,\ell+1}){D}_{\ell,\ell}\big{)}{W}_{\ell}^\dagger {H}_{\ell},\\[1mm]
T_{k,N}=0,\\[1mm]
\ell\in \mathbb{T}_{k},
\end{array}
\right.\\[1mm]
{W}_k{W}_k^\dagger{H}_k-{H}_k=0,\\[1mm]
k\in \mathbb{T},
\end{array}
\right.
\end{eqnarray}
and
\begin{eqnarray}\label{pi-2-s2}
\left\{
\begin{array}{l}
\left\{
\begin{array}{l}
\pi_{k,\ell}={A}^T_{k,\ell}{P}_{k,\ell+1}\big{(}f_{k,\ell}-{B}_{k,\ell}{W}_{\ell}^\dagger \beta_{\ell}\big{)}+{A}^T_{k,\ell}{T}_{k,\ell+1}\big{(}f_{\ell,\ell}-{B}_{\ell,\ell} {W}_{\ell}^\dagger \beta_{\ell}\big{)}\\[1mm]
\hphantom{\pi_{k,\ell}=}+{C}^T_{k,\ell}P_{k,\ell+1}\big{(}d_{k,\ell}-{D}_{k,\ell}{W}_{\ell}^\dagger \beta_{\ell}\big{)}+{C}^T_{k,\ell}T_{k,\ell+1}\big{(}d_{\ell,\ell}-{D}_{\ell,\ell}{W}_{\ell}^\dagger \beta_{\ell}\big{)}\\[1mm]
\hphantom{\pi_{k,\ell}=}+{A}^T_{k,\ell}\pi_{k,\ell+1}+q_{k,\ell},\\[1mm]
\pi_{k,N}=g_k,
\end{array}
\right.\\[1mm]
{W}_k{W}_k^\dagger\beta_k-\beta_k=0,\\[1mm]
k\in \mathbb{T},
\end{array}
\right.
\end{eqnarray}
where
\begin{eqnarray*}\label{feedback-gain-k-in}
\left\{
\begin{array}{l}
{W}_{k}={R}_{k,k}+{B}_{k,k}^T\big{(}{P}_{k,k+1}+{T}_{k,k+1}\big{)}{B}_{k,k}+{D}_{k,k}^T\big{(}{P}_{k,k+1}+T_{k,k+1}\big{)}{D}_{k,k},\\[1mm]
{H}_{k}={B}_{k,k}^T\big{(}{P}_{k,k+1}+{T}_{k,k+1}\big{)}{A}_{k,k}+{D}_{k,k}^T\big{(}{P}_{k,k+1}+T_{k,k+1}\big{)}{C}_{k,k},\\[1mm]
\beta_{k}={B}_{k,k}^T\big{[}\big{(}{{P}}_{k,k+1}+{T}_{k,k+1}\big{)}f_{k,k}+\pi_{k,k+1}\big{]}+{D}_{k,k}^T\big{(}{P}_{k,k+1}+T_{k,k+1}\big{)}d_{k,k}+\rho_{k,k}, \\[1mm] k\in \mathbb{T}.
\end{array}
\right.
\end{eqnarray*}

\begin{corollary}\label{corollary-Necessary-sufficient-2}
For any $t\in \mathbb{T}$ and any $x\in L^2_\mathcal{F}(t;
\mathbb{R}^n)$, Problem (LQ)$_{s2}$ for the initial pair $(t,x)$
admits an open-loop equilibrium pair if and only if (\ref{P-2-s2}), (\ref{T-2-s2}) and (\ref{pi-2-s2}) are solvable.
\end{corollary}

In Theorem 2.2 of \cite{Li-Ni-Zhang}, a necessary and sufficient
condition to the existence of open-loop equilibrium pair is
presented for the following system
\begin{eqnarray}\label{system-1-in-2}
\left\{\begin{array}{l}
X_{k+1}=\big{(}A_{k}X_k+B_{k}u_k+f_{k}\big{)}+\big{(}C_{k}X_k+D_{k}u_k+d_{k}\big{)}w_k, \\[1mm]
X_t=x,~~k\in \mathbb{T}_t,~~t\in \mathbb{T}.
\end{array}
\right.
\end{eqnarray}
The matrices in (\ref{system-1-in-2}) are independent of the initial time. Hence, Corollary \ref{corollary-Necessary-sufficient-2} is an
extension of Theorem 2.2 in \cite{Li-Ni-Zhang}.

\section{Example}\label{section-4}

In this section, we will use an example to illustrate the theory on
solving Problem (LQ).

\textbf{Example 4.1} Consider a version of Problem (LQ) with $N=2$, whose system matrices
and weighting matrices are given blow
\begin{eqnarray*}
&&\hspace{-1.3em}A_{0,0}=\left[
\begin{array}{cc}
3.3 & 0.41\\
-1.3& 1.9
\end{array}
\right],~~
A_{0,1}=\left[
\begin{array}{cc}
5.12 & -0.35\\
1.31& 2.03
\end{array}
\right],~~\bar{A}_{0,0}=\left[
\begin{array}{cc}
3.34 & -1.01\\
1.43& 2.03
\end{array}
\right],\\[1mm]
&&\hspace{-1.3em}\bar{A}_{0,1}=\left[
\begin{array}{cc}
3.45 & -0.3\\
1.2& 4
\end{array}
\right],
B_{0,0}=\left[
\begin{array}{cc}
3.5 & 1.6\\
-0.2& 3
\end{array}
\right],~~B_{0,1}=\left[
\begin{array}{cc}
4.45 & 2.36\\
-1.2& 5
\end{array}
\right],~\\[1mm]
&&\hspace{-1.3em}\bar{B}_{0,0}=\left[
\begin{array}{cc}
3.2 & 0.32\\
1.5&  3
\end{array}
\right],~~\bar{B}_{0,1}=\left[
\begin{array}{cc}
3.65 & -0.3\\
-0.42& 5.6
\end{array}
\right],
C_{0,0}=\left[
\begin{array}{cc}
5.6 & 1\\
0.73&  7.8
\end{array}
\right],\\[1mm]
&&\hspace{-1.3em}C_{0,1}=\left[
\begin{array}{cc}
5 & 0.73\\
-0.47& 5.2
\end{array}
\right],~~\bar{C}_{0,0}=\left[
\begin{array}{cc}
5.6 & 1\\
0.73&  7.8
\end{array}
\right],~~\bar{C}_{0,1}=\left[
\begin{array}{cc}
5 & 0.73\\
-0.47& 5.2
\end{array}
\right],\\[1mm]
&&\hspace{-1.3em}D_{0,0}=\left[
\begin{array}{cc}
6 & 1.63\\
-1.37& 7
\end{array}
\right],~~
D_{0,1}=\left[
\begin{array}{cc}
4 & 0.93\\
1.07& 3
\end{array}
\right],~~\bar{D}_{0,0}=\left[
\begin{array}{cc}
4.6 & 0.63\\
-1.57& 6.4
\end{array}
\right],\\[1mm]
&&\hspace{-1.3em}
\bar{D}_{0,1}=\left[
\begin{array}{cc}
4.4 & 1.93\\
2.34& 5.63
\end{array}
\right],~~
A_{1,1}=\left[
\begin{array}{cc}
8.5 & 3.03\\
-2.23& 7.2
\end{array}
\right],~~
\bar{A}_{1,1}=\left[
\begin{array}{cc}
5.67 & 1.93\\
-1.16 & 6.54
\end{array}
\right],\\[1mm]
&&\hspace{-1.3em}
B_{1,1}=\left[
\begin{array}{cc}
7.35 & -2.35\\
-3.38& 6.32
\end{array}
\right],~~
\bar{B}_{1,1}=\left[
\begin{array}{cc}
5.67 & 1.93\\
-1.16 & 6.54
\end{array}
\right],
C_{1,1}=\left[
\begin{array}{cc}
2.5 & 3.03\\
-4.23& 6.2
\end{array}
\right],~~\\[1mm]
&&\hspace{-1.3em}
\bar{C}_{1,1}=\left[
\begin{array}{cc}
10.17 & 5.93\\
-6.16 & 7.54
\end{array}
\right],~~D_{1,1}=\left[
\begin{array}{cc}
8.56 & -4.75\\
-2.8& 7
\end{array}
\right],~~
\bar{D}_{1,1}=\left[
\begin{array}{cc}
-8.72 & 2.43\\
1.16 & -6.54
\end{array}
\right],\\[1mm]
&&\hspace{-1.3em}Q_{0,0}=\left[
\begin{array}{cc}
-1 & 0.8\\
0.8&  -1.6
\end{array}
\right],~~Q_{0,1}=\left[
\begin{array}{cc}
4 & 0\\
0& 0
\end{array}
\right],~~\bar{Q}_{0,0}=\left[
\begin{array}{cc}
-0.5 & -0.1\\
-0.1&  1
\end{array}
\right],~~\bar{Q}_{0,1}=\left[
\begin{array}{cc}
-2 & 0\\
0& -3
\end{array}
\right],\\[1mm]
&&\hspace{-1.3em}R_{0,0}=\left[
\begin{array}{cc}
-0.5 & 0\\
0&  1
\end{array}
\right],~~R_{0,1}=\left[
\begin{array}{cc}
1 & 0\\
0&  -2
\end{array}
\right],~~\bar{R}_{0,0}=\left[
\begin{array}{cc}
0 & 0\\
0&  0
\end{array}
\right],~~\bar{R}_{0,1}=\left[
\begin{array}{cc}
-2 & 0\\
0&  2
\end{array}
\right],\\[1mm]
&&\hspace{-1.3em}Q_{1,1}=\left[
\begin{array}{cc}
2 & 0.1\\
0.1& 5
\end{array}
\right],~~
\bar{Q}_{1,1}=\left[
\begin{array}{cc}
-1 & 0.1\\
0.1& -3
\end{array}
\right],~~R_{1,1}=\left[
\begin{array}{cc}
4 & -0.3\\
-0.3& -2
\end{array}
\right],~~
\bar{R}_{1,1}=\left[
\begin{array}{cc}
-7 & -1.3\\
-1.3& -4
\end{array}
\right],
\\[1mm]
&&\hspace{-1.3em}G_0=\left[
\begin{array}{cc}
1 & 0\\
0&  2
\end{array}
\right],~~G_1=\left[
\begin{array}{cc}
2 & -0.3\\
-0.3&  3
\end{array}
\right],~~\bar{G}_0=\left[
\begin{array}{cc}
2 & 0\\
0&  3
\end{array}
\right],~~\bar{G}_1=\left[
\begin{array}{cc}
-0.5 & -0.2\\
-0.2&  1
\end{array}
\right],\\[1mm]
&&\hspace{-1.3em}f_{0,0}=\left[
\begin{array}{cc}
-0.5 \\ -1
\end{array}
\right],~~f_{0,1}=\left[
\begin{array}{cc}
-1.34 \\ 2.5
\end{array}
\right],~~d_{0,0}=\left[
\begin{array}{cc}
1.32 \\ 2.79
\end{array}
\right],~~d_{0,1}=\left[
\begin{array}{cc}
-0.35 \\ 8.9
\end{array}
\right],~~\\[1mm]
&&\hspace{-1.3em}f_{1,1}=\left[
\begin{array}{cc}
1 \\ 2
\end{array}
\right],~~
~~d_{1,1}=\left[
\begin{array}{cc}
0 \\ 1
\end{array}
\right],~~q_{0,0}=\left[
\begin{array}{cc}
-0.85 \\ -1.8
\end{array}
\right],~~
q_{0,1}=\left[
\begin{array}{cc}
2 \\ 7
\end{array}
\right],~~\rho_{0,0}=\left[
\begin{array}{cc}
3.2 \\ 2.1
\end{array}
\right],~~~~\\[1mm]
&&\hspace{-1.3em}\rho_{0,1}=\left[
\begin{array}{cc}
1.42 \\ 2.71
\end{array}
\right],~~q_{1,1}=\left[
\begin{array}{cc}
6 \\ 8
\end{array}
\right],~~\rho_{1,1}=\left[
\begin{array}{cc}
6.2 \\ -5.7
\end{array}
\right],~~g_{0}=\left[
\begin{array}{cc}
5.6 \\ 7.8
\end{array}
\right],~~g_{1}=\left[
\begin{array}{cc}
-9 \\ 8.7
\end{array}
\right].
\end{eqnarray*}
Note that $Q_{k,\ell}, \mathcal{Q}_{k,\ell}$, $R_{k,\ell},
\mathcal{R}_{k,\ell}$, $k=0,1, \ell=k,1$, are not fully nonnegative
definite, since for example $Q_{0,0}$ is negative definite and $R_{0,0}$ is
indefinite.

By the iterations of (\ref{p-3}), (\ref{T-3}) and (\ref{pi-3}), we
can get the values of the solutions with
\begin{eqnarray*}
&&\mathcal{R}_{1,1}+\mathcal{B}_{1,1}^T\mathcal{P}_{1,2}\mathcal{B}_{1,1}
+\mathcal{D}_{1,1}^TP_{1,2}\mathcal{D}_{1,1}=\left[
\begin{array}{cc}
  400.8004~&~ -330.6524\\
 -330.6524~&~  673.2241
\end{array}
\right]>0,\\[1mm]
&&\mathcal{R}_{0,0}+\mathcal{B}_{0,0}^T\mathcal{P}_{0,1}\mathcal{B}_{0,0}
+\mathcal{D}_{0,0}^TP_{0,1}\mathcal{D}_{0,0}=\left[
\begin{array}{cc}
    24209~&~    11560\\
    11560 ~&~   28652
\end{array}
\right]>0,\\[1mm]
&&\mathcal{W}_{1}=\left[
\begin{array}{cc}
  400.8004~&~ -330.6524\\
 -330.6524 ~&~ 673.2241
\end{array}
\right],~~\mathcal{W}_{0}=\left[
\begin{array}{cc}
    12637 ~&~   932\\
   -6334 ~&~   3464
\end{array}
\right].
\end{eqnarray*}
The sets of eigenvalues of $\mathcal{W}_{1}$ an  $\mathcal{W}_0$ are $\{179.4026, 894.6219\}$ and $\{11940, 4160\}$, respectively. Hence, $\mathcal{W}_1$ and $\mathcal{W}_0$ are both invertible. Therefore, the corresponding (\ref{p-3}), (\ref{T-3}) and (\ref{pi-3})
are solvable, and for any initial pair $(t,x)$ with $t=0,1, x\in L^2_\mathcal{F}(t; \mathbb{R}^2)$
the considered LQ problem admits an open-loop equilibrium pair.
Furthermore, an open-loop equilibrium control for the initial pair $(0,x)$ is given by
\begin{eqnarray*}
u^{0,x,*}_k=-\mathcal{W}_k^\dagger \mathcal{H}_kX^{0,x,*}_k-\mathcal{W}_k^\dagger \beta_k, ~~k=0,1,
\end{eqnarray*}
where
\begin{eqnarray*}
&&\mathcal{W}_{1}^\dagger \mathcal{H}_1=\left[
\begin{array}{cc}
    1.1320~  &  ~0.1179\\
    0.0254~  &  ~1.0388
\end{array}
\right],~~\mathcal{W}_{0}^\dagger \mathcal{H}_0=\left[
\begin{array}{cc}
    0.8661 ~&~  -0.4704\\
    0.0520 ~&~   0.9824
\end{array}
\right],\\[1mm]
&&\mathcal{W}_{1}^\dagger \beta_1=\left[
\begin{array}{cc}
   -0.3381\\
    0.1433
\end{array}
\right],~~\mathcal{W}_{0}^\dagger \beta_0=\left[
\begin{array}{cc}
   -0.2003\\
   -0.1582
\end{array}
\right],
\end{eqnarray*}
and
\begin{eqnarray*}
\left\{
\begin{array}{l}
X^{0,x,*}_{k+1} = \big{[}\mathcal{A}_{k,k}X^{0,x,*}_{k}+\mathcal{B}_{k,k}u^{0,x,*}_k+f_{k,k}\big{]}\\[1mm]
\hphantom{X^{0,x,*}_{k+1} = }+\big{[}\mathcal{C}_{k,k}X^{0,x,*}_{k} +\mathcal{D}_{k,k}u^{0,x,*}_k+d_{k,k}\big{]}w_k,\\[1mm]
X^{0,x,*}_{0} = x,~~ k \in  \{0,1\}.
\end{array}
\right.
\end{eqnarray*}

\section{Conclusion}\label{section-5}

In this paper, the open-loop time-consistent equilibrium control is
investigated for a kind of mean-field stochastic LQ problem, where both the system matrices
and the weighting matrices are depending on the initial time, and the conditional
expectations of the control and state enter quadratically into the cost functional.
Necessary and sufficient conditions are presented for both the case with
a fixed initial pair and the case with all the initial pairs. Furthermore, a set of constrained GDREs and two sets of constrained LDEs are introduced to characterize the closed-loop form
of open-loop equilibrium control. Note that this paper is concerned with the time-consistency of open-loop control. For future research, the time-consistency of the strategy should be studied.



\section*{Appendix}

\subsection*{A. Proof of Lemma \ref{Lemma-difference}}\label{appendix-A}         

\emph{Proof}. Let us replace $u_k$ with $u_k+\lambda\bar{u}_k$ in
the forward S$\Delta$E of (\ref{X-Z-1}), and denote its solution by
${X}^{k,\lambda}$. Then, we have
\begin{eqnarray*}
\left\{
\begin{array}{l}
\frac{{X}^{k,\lambda}_{\ell+1}-X^{k,u_k}_{\ell+1}}{\lambda}=\Big{(}A_{k,\ell}\frac{{X}_\ell^{k,\lambda}-X^{k,u_k}_\ell}{\lambda}+\bar{A}_{k,\ell}\frac{\mathbb{E}_k{X}_\ell^{k,\lambda}-\mathbb{E}_kX^{k,u_k}_\ell}{\lambda}\Big{)}\\[1mm]
\hphantom{\frac{{X}^{k,\lambda}_{\ell+1}-X^{k,u_k}_{\ell+1}}{\lambda}=} +\Big{(}C_{k,\ell}\frac{{X}_\ell^{k,\lambda}-X^{k,u_k}_\ell}{\lambda}+\bar{C}_{k,\ell}\frac{\mathbb{E}_k{X}_\ell^{k,\lambda}-\mathbb{E}_kX^{k,u_k}_\ell}{\lambda}\Big{)}w_\ell,\\[1mm]
\frac{{X}^{k,\lambda}_{k+1}-X^{k,u_k}_{k+1}}{\lambda}=\Big{(}A_{k,k}\frac{{X}_k^{k,\lambda}-X^{k,u_k}_k}{\lambda}+ \bar{A}_{k,k}\frac{\mathbb{E}_k{X}_k^{k,\lambda}-\mathbb{E}_kX^{k,u_k}_k}{\lambda}+B_{k,k}\bar{u}_k+\bar{B}_{k,k}\bar{u}_k\Big{)}\\[1mm]
\hphantom{\frac{{X}^{k,\lambda}_{\ell+1}-X^{k,u_k}_{\ell+1}}{\lambda}=}+\Big{(}C_{k,k}\frac{{X}_k^{k,\lambda}-X^{k,u_k}_k}{\lambda}+\bar{C}_{k,k}\frac{\mathbb{E}_k{X}_k^{k,\lambda}-\mathbb{E}_kX^{k,u_k}_k}{\lambda}+D_{k,k}\bar{u}_k+\bar{D}_{k,k}\mathbb{E}_k\bar{u}_k\Big{)}w_k,\\
\frac{{X}_k^{k,\lambda}-X^{k,u_k}_k}{\lambda}=0,~~~\ell\in \mathbb{T}_{k+1}.
\end{array}
\right.
\end{eqnarray*}
Denoting $\frac{{X}_\ell^{k,\lambda}-X^{k,u_k}_\ell}{\lambda}$ by ${Y}^{k,\bar{u}_k}_\ell$, we get
\begin{eqnarray}\label{system-y}
\left\{
\begin{array}{l}
{Y}^{k,\bar{u}_k}_{\ell+1}=A_{k,\ell}{Y}^{k,\bar{u}_k}_\ell+\bar{A}_{k,\ell}\mathbb{E}_k{Y}^{k,\bar{u}_k}_\ell+\big{(}C_{k,\ell}{Y}^{k,\bar{u}_k}_\ell+\bar{C}_{k,\ell}\mathbb{E}_k{Y}^{k,\bar{u}_k}_\ell\big{)}w_\ell,\\[1mm]
Y^{k,\bar{u}_k}_{k+1}=(B_{k,k}+\bar{B}_{k,k})\bar{u}_k+{(}D_{k,k}+\bar{D}_{k,k}{)}\bar{u}_k w_k,\\
{Y}^{k,\bar{u}_k}_k=0,~~\ell\in \mathbb{T}_{k+1}.
\end{array}
\right.
\end{eqnarray}
Here, we have used the fact $\mathbb{E}_ku_k=u_k$. Note that $X^{k, \lambda}_\ell=X^{k,u}_\ell
+\lambda Y^{k,\bar{u}_k}_\ell$, $ \forall \ell \in \mathbb{T}_k$. Then, we have
\begin{eqnarray}\label{appendix-A-J}
&&\hspace{-2em}{J}(k,\zeta; (u_k+\lambda \bar{u}_k,u|_{\mathbb{T}_{k+1}}))-{J}(k,\zeta; u)\nonumber\\[1mm]
&&\hspace{-2em} =\sum_{\ell=k}^{N-1}\dbE_k\Big{\{}(X^{k,u_k}_\ell+\lambda Y^{k,\bar{u}_k}_\ell)^TQ_{k,\ell}(X^{k,u_k}_\ell+\lambda Y^{k,\bar{u}_k}_\ell)+[\mathbb{E}_k(X^{k,u_k}_\ell+\lambda Y^{k,\bar{u}_k}_\ell)]^T \bar{Q}_{k,\ell}\mathbb{E}_k(X^{k,u_k}_\ell+\lambda Y^{k,\bar{u}_k}_\ell)\nonumber \\
&&\hspace{-2em}\hphantom{=}+2q_{k,\ell}^T (X^{k,u_k}_\ell+\lambda Y^{k,\bar{u}_k}_\ell)-(X^{k,u_k}_\ell)^TQ_{k,\ell}X^{k,u_k}_\ell-[\mathbb{E}_kX^{k,u_k}_\ell]^T \bar{Q}_{k,\ell}\mathbb{E}_\ell X^{k,u_k}_\ell-2q_{k,\ell}^T X^{k,u_k}_\ell\Big{\}}\nonumber\\[1mm]
&&\hspace{-2em}\hphantom{=}+(u_k+\lambda \bar{u}_k)^T(R_{k,k}+\bar{R}_{k,k})(u_k+\lambda \bar{u}_k)+
2\rho_{k,\ell}^T (u_k+\lambda \bar{u}_k)-u_k^T(R_{k,k}+\bar{R}_{k,k})u_k-2\rho_{k,\ell}^T u_k\nonumber\\[1mm]
&&\hspace{-2em}\hphantom{=}+[\mathbb{E}_k(X^{k,u_k}_N+\lambda Y_N^{k,\bar{u}_k})]^T \bar{G}_k\mathbb{E}_k(X^{k,u_k}_{N}+\lambda Y_N^{k,\bar{u}_k})+\mathbb{E}_k\big{[}(X_N^{k,u_k}+\lambda Y_N^{k,\bar{u}_k})^TG_k(X^{k,u_k}_N+\lambda Y_N^{k,\bar{u}_k})\big{]}\nonumber\\[1mm]
&&\hspace{-2em}\hphantom{=}+2\mathbb{E}_k[g_k^T (X^{k,u_k}_N+\lambda Y_N^{k,\bar{u}_k})] -\mathbb{E}_k\big{[}(X_N^{k,u_k})^TG_kX^{k,u_k}_N\big{]}-(\mathbb{E}_kX^{k,u_k}_N)^T \bar{G}_k\mathbb{E}_kX^{k,u_k}_{N}-2\mathbb{E}_kg_k^T X^{k,u_k}_N\nonumber\\[1mm]
&&\hspace{-2em} =2\lambda \Big{\{} \sum_{\ell=k}^{N-1}\dbE_k\Big{[}(X^{k,u_k}_\ell )^TQ_{k,\ell} Y^{k,\bar{u}_k}_\ell+q_{k,\ell}^T Y^k_\ell+[\mathbb{E}_k X^{k,u_k}_\ell]^T \bar{Q}_{k,\ell}\mathbb{E}_k Y^{k,\bar{u}_k}_\ell\Big{]}+u_k^T(R_{k,k}+\bar{R}_{k,k})\bar{u}_k\nonumber\\[1mm]
&&\hspace{-2em}\hphantom{=}+\rho_{k,\ell}^T \bar{u}_k+\mathbb{E}_k\big{[}(X_N^{k,u_k})^TG_k Y_N^{k,\bar{u}_k}\big{]}+\mathbb{E}_k[g_k^T  Y_N^{k,\bar{u}_k}]+[\mathbb{E}_kX^{k,u_k}_N]^T \bar{G}_k\mathbb{E}_k Y_N^{k,\bar{u}_k}\Big{\}}\nonumber\\
&&\hspace{-2em}\hphantom{=}+ \lambda^2\Big{\{} \sum_{\ell=k}^{N-1}\mathbb{E}_k\Big{[}(Y^{k,\bar{u}_k}_\ell)^TQ_{k,\ell} Y_\ell^{k,\bar{u}_k}+(\mathbb{E}_kY^{k,\bar{u}_k}_\ell)^T\bar{Q}_{k,\ell} \mathbb{E}_kY_\ell^{k,\bar{u}_k}\Big{]}+\mathbb{E}_k\big{[} \bar{u}_k^T(R_{k,k}+\bar{R}_{k,k})\bar{u}_k \big{]}\nonumber\\
&&\hspace{-2em}\hphantom{=} +\mathbb{E}_k\big{[}(Y_N^{k,\bar{u}_k})^T G_{k} {Y}_N^{k,\bar{u}_k}\big{]}
+(\mathbb{E}_kY_N^{k,\bar{u}_k})^T \bar{G}_{k} \mathbb{E}_k{Y}_N^{{k,\bar{u}_k}}\Big{\}}.
\end{eqnarray}
On the other hand, we have
\begin{eqnarray*}
&&\hspace{-2em}\sum_{\ell=k}^{N-1}\dbE_k\Big{[}(X^{k,u_k}_\ell )^TQ_{k,\ell} Y^{k,\bar{u}_k}_\ell+q_{k,\ell}^T Y^{k,\bar{u}_k}_\ell+[\mathbb{E}_k X^{k,u_k}_\ell]^T \bar{Q}_{k,\ell}\mathbb{E}_k Y^{k,\bar{u}_k}_\ell \Big{]}+u_k^T(R_{k,k}+\bar{R}_{k,k})\bar{u}_k\nonumber\\[1mm]
&&\hspace{-2em}+\rho_{k,\ell}^T \bar{u}_k+\mathbb{E}_k\big{[}(X_N^{k,u_k})^TG_k Y_N^{k,\bar{u}_k}\big{]}+[\mathbb{E}_kX^{k,u_k}_N]^T \bar{G}_k\mathbb{E}_k Y_N^{k,\bar{u}_k}+\mathbb{E}_k[g_k^T  Y_N^{k,\bar{u}_k}]\nonumber\\
%
%
%
%
%
&&\hspace{-2em}=\sum_{\ell=k}^{N-1}\mathbb{E}_k\Big{[}\Big{(}Q_{k,\ell}(X_\ell^{k,u_k}-\mathbb{E}_kX_\ell^{k,u_k})+A_{k,\ell}^T(\mathbb{E}_\ell Z_{\ell+1}^{k,u_k}-\mathbb{E}_kZ^{k,u_k}_{\ell+1}) \\[1mm] &&\hspace{-2em}\hphantom{=}+C_{k,\ell}^T\big{(}\mathbb{E}_\ell(Z_{\ell+1}^{k,u_k}w_\ell)-\mathbb{E}_k(Z_{\ell+1}^{k,u_k}w_\ell)\big{)}-(Z_{\ell}^{k,u_k}-\mathbb{E}_kZ_\ell^{k,u_k})\Big{)}^T (Y^{k,\bar{u}_k}_\ell-\mathbb{E}_kY^{k,\bar{u}_k}_\ell)\\[1mm]
&&\hspace{-2em}\hphantom{=}+\Big{(}(Q_{k,\ell}+\bar{Q}_{k,\ell})\mathbb{E}_kX_\ell^{k,u_k}+q_{k,\ell}+(A_{k,\ell}+\bar{A}_{k,\ell})^T\mathbb{E}_kZ_{\ell+1}^{k,u_k} \nonumber \\[1mm]
&&\hspace{-2em} \hphantom{=} +(C_{k,\ell}+\bar{C}_{k,\ell})^T\mathbb{E}_k(Z_{\ell+1}^{k,u_k}w_\ell)-\mathbb{E}_kZ_\ell^{k,u_k}\Big{)}^T\mathbb{E}_kY_\ell^{k,\bar{u}_k}\Big{]}\nonumber \\[1mm]
&&\hspace{-2em}\hphantom{=}+\Big{[}(R_{k,k}+\bar{R}_{k,k})u_k+(B_{k,k}+\bar{B}_{k,k})^T\mathbb{E}_kZ^{k,u_k}_{k+1}+(D_{k,k}+\bar{D}_{k,k})^T\mathbb{E}_k(Z_{k+1}^{k,u_k}w_k)+\rho_{k,k}\Big{]}^T\bar{u}_k\nonumber\\[1mm]
&&\hspace{-2em}=\Big{[}(R_{k,k}+\bar{R}_{k,k})u_k+(B_{k,k}+\bar{B}_{k,k})^T\mathbb{E}_kZ^{k,u_k}_{k+1}+(D_{k,k}+\bar{D}_{k,k})^T\mathbb{E}_k(Z_{k+1}^{k,u_k}w_k)+\rho_{k,k}\Big{]}^T\bar{u}_k.
\end{eqnarray*}
This together with (\ref{appendix-A-J}) implies the conclusion.\endpf

\subsection*{B. Proof of Theorem \ref{Theorem-Equivalentce-open-loop}}\label{appendix-B}         

\emph{Proof.} (i)$\Rightarrow$(ii). Let $(X^{t,x,*}, u^{t,x,*})$ be
an equilibrium pair. As (\ref{system-adjoint}) is a decoupled FBS$\Delta$E,
(\ref{system-adjoint}) is solvable. From  (\ref{appendix-A-J-0}) we
have
\begin{eqnarray}\label{appendix-B-1}
&&\hspace{-2em}{J}(k,X^{t,x,*}_k; (u^{t,x,*}_k+\lambda \bar{u}_k, u^{t,x,*}|_{\mathbb{T}_{k+1}}))-{J}(k,X^{t,x,*}_k; u^{t,x,*})\nonumber\\
&&\hspace{-2em} =2\lambda \Big{[}(R_{k,k}+\bar{R}_{k,k})u_k^{t,x,*}+(B_{k,k}+\bar{B}_{k,k})^T\mathbb{E}_kZ^{k,t,x}_{k+1}+(D_{k,k}+\bar{D}_{k,k})^T\mathbb{E}_k(Z_{k+1}^{k,t,x}w_k)+\rho_{k,k}\Big{]}^T\bar{u}_k.\nonumber\\[1mm]
&&\hspace{-2em}\hphantom{=}+ \lambda^2\Big{\{} \sum_{\ell=k}^{N-1}\mathbb{E}_k\Big{[}(Y^{k,\bar{u}_k}_\ell)^TQ_{k,\ell} Y_\ell^{k,\bar{u}_k}+(\mathbb{E}_kY^{k,\bar{u}_k}_\ell)^T\bar{Q}_{k,\ell} \mathbb{E}_kY_\ell^{k,\bar{u}_k}\Big{]}+\mathbb{E}_k\big{[} \bar{u}_k^T(R_{k,k}+\bar{R}_{k,k})\bar{u}_k \big{]}\nonumber\\
&&\hspace{-2em}\hphantom{=} +\mathbb{E}_k\big{[}(Y_N^{k,\bar{u}_k})^T G_{k} {Y}_N^{{k,\bar{u}_k}}\big{]}+(\mathbb{E}_kY_N^{k,\bar{u}_k})^T \bar{G}_{k} \mathbb{E}_k{Y}_N^{{k,\bar{u}_k}}\Big{\}}\nonumber \\
&&\hspace{-2em}\geq 0.
\end{eqnarray}
Noting that (\ref{appendix-B-1}) holds for any $\lambda \in
\mathbb{R}$ and $\bar{u}_k\in L^2_{\mathcal{F}}(k; \mathbb{R}^m)$,
we have (\ref{stationary-condition}) and (\ref{convex}). In fact, if
(\ref{convex}) was not satisfied, then there would be a $\bar{u}_k$
such that $\ds\lim_{\lambda\mapsto \infty}\bar{J}(k,X^{t,x,*}_k;
u^{t,x,*}_k+\lambda \bar{u}_k)-\bar{J}(k,X^{t,x,*}_k;
u^{t,x,*}_k)=-\infty$. This is impossible. Furthermore, for any given $\bar{u}_k\in L^2_{\mathcal{F}}(k; \mathbb{R}^m)$, denote
\begin{eqnarray*}
&&\hspace{-2.5em}\delta_1 \triangleq \Big{[}(R_{k,k}+\bar{R}_{k,k})u_k^{t,x,*}+(B_{k,k}+\bar{B}_{k,k})^T\mathbb{E}_kZ^{k,t,x}_{k+1}\\[1mm]
&&\hspace{-2.5em}\hphantom{\delta_1 \triangleq} +(D_{k,k}+\bar{D}_{k,k})^T\mathbb{E}_k(Z_{k+1}^{k,t,x}w_k)+\rho_{k,k}\Big{]}^T\bar{u}_k\neq 0,
\end{eqnarray*}
and
\begin{eqnarray*}
&&\hspace{-2em}\delta_2=\sum_{\ell=k}^{N-1}\mathbb{E}_k\Big{[}(Y^{k,\bar{u}_k}_\ell)^TQ_{k,\ell} Y_\ell^{k,\bar{u}_k}+(\mathbb{E}_kY^{k,\bar{u}_k}_\ell)^T\bar{Q}_{k,\ell} \mathbb{E}_kY_\ell^{k,\bar{u}_k}\Big{]}\\[1mm]
&&\hspace{-2em}\hphantom{\delta_2=}+\mathbb{E}_k\big{[} \bar{u}_k^T(R_{k,k}+\bar{R}_{k,k})\bar{u}_k \big{]}+\mathbb{E}_k\big{[}(Y_N^{k,\bar{u}_k})^T G_{k} {Y}_N^{{k,\bar{u}_k}}\big{]}+(\mathbb{E}_kY_N^{k,\bar{u}_k})^T \bar{G}_{k} \mathbb{E}_k{Y}_N^{{k,\bar{u}_k}}.
\end{eqnarray*}
When $\delta_2=0$, we select $\lambda=-\delta_1$, which together with
(\ref{appendix-B-1}) implies
\begin{eqnarray*}
{J}(k,X^{t,x,*}_k; (u^{t,x,*}_k+\lambda \bar{u}_k, u^{t,x,*}|_{\mathbb{T}_{k+1}}))-{J}(k,X^{t,x,*}_k; u^{t,x,*})=-\delta_1^2<0.
\end{eqnarray*}
This is impossible. When  $\delta_2\neq 0$ (which is positive), we
select $\lambda =\theta \delta_1<0$  with
$\theta=-\frac{1}{\delta_2} $. In this case, we have
\begin{eqnarray*}
{J}(k,X^{t,x,*}_k; (u^{t,x,*}_k+\lambda \bar{u}_k, u^{t,x,*}|_{\mathbb{T}_{k+1}}))-{J}(k,X^{t,x,*}_k; u^{t,x,*})=2\theta \delta_1^2+\theta^2 \delta_1^2\delta_2=\delta_1^2\theta<0,
\end{eqnarray*}
which contradicts (\ref{appendix-B-1}).

(ii)$\Rightarrow$(i). In this case, for any $\lambda \in \mathbb{R}$
and $\bar{u}_k\in L^2_{\mathcal{F}}(k; \mathbb{R}^m)$ we have
\begin{eqnarray}\label{appendix-B-2}
&&\hspace{-2em}{J}(k,X^{t,x,*}_k; (u^{t,x,*}_k+\lambda \bar{u}_k, u^{t,x,*}|_{\mathbb{T}_{k+1}}))-{J}(k,X^{t,x,*}_k; u^{t,x,*})\nonumber\\
&&\hspace{-2em} =\lambda^2\Big{\{} \sum_{\ell=k}^{N-1}\mathbb{E}_k\Big{[}(Y^{k,\bar{u}_k}_\ell)^TQ_{k,\ell} Y_\ell^{k,\bar{u}_k}+(\mathbb{E}_kY^{k,\bar{u}_k}_\ell)^T\bar{Q}_{k,\ell} \mathbb{E}_kY_\ell^{k,\bar{u}_k}\Big{]}\nonumber\\
&&\hspace{-2em}\hphantom{=}+\mathbb{E}_k\big{[} \bar{u}_k^T(R_{k,k}+\bar{R}_{k,k})\bar{u}_k \big{]}+\mathbb{E}_k\big{[}(Y_N^{k,\bar{u}_k})^T G_{k} {Y}_N^{{k,\bar{u}_k}}\big{]}+(\mathbb{E}_kY_N^{k,\bar{u}_k})^T \bar{G}_{k} \mathbb{E}_k{Y}_N^{{k,\bar{u}_k}}\Big{\}}\nonumber \\
&&\hspace{-2em}\geq 0.
\end{eqnarray}
For $k=t$, from (\ref{system-1}) we have
\begin{eqnarray}\label{system-1-t=0}
\left\{\begin{array}{l}
X^{t}_{t+1}=\big{[}(A_{t,t}+\bar{A}_{t,t})X^t_t+(B_{t,t}+\bar{B}_{t,t})u^{t,x,*}_t+f_{t,t}\big{]}\\[1mm]
\hphantom{X^{t}_{t+1}=}+\big{[}(C_{t,t}+\bar{C}_{t,t})X^t_t+(D_{t,t}+\bar{D}_{t,t})u^{t,x,*}_t+d_{t,t}\big{)}w_t, \\[1mm]
X^{t,x,*}_{t} = x,
\end{array}\right.
\end{eqnarray}
and for any $u_t\in L^2_\mathcal{F}(t; \mathbb{R}^m)$,
\begin{eqnarray}\label{open-loop-equilibrium-t=0}
J(t,x; u^{t,x,*})\leq J(t, x; (u_t, u^{t,x,*}|_{\mathbb{T}_{t+1}})).
\end{eqnarray}
We now move to the case of $k=t+1$. In this case, the starting point of the
state is $X_{t+1}^{t,x,*}$. Hence, we have
\begin{eqnarray}\label{system-1-t=1}
\left\{\begin{array}{l}
X^{t+1}_{t+2}=\big{[}(A_{t+1,t+1}+\bar{A}_{t+1,t+1})X^{t+1}_{t+1}+(B_{t+1,t+1}+\bar{B}_{t+1,t+1})u^{t,x,*}_{t+1}+f_{t+1,t+1}\big{]}\\[1mm]
\hphantom{X^{t+1}_{t+2}=}+\big{[}(C_{t+1,t+1}+\bar{C}_{t+1,t+1})X^{t+1}_{t+1}+(D_{t+1,t+1}+\bar{D}_{t+1,t+1})u^{t,x,*}_{t+1}+d_{t+1,t+1}\big{)}w_{t+1}, \\[1mm]
X^{t+1}_{t+1} = X_{t+1}^{t},
\end{array}\right.
\end{eqnarray}
and for any $u_{t+1}\in L^2_\mathcal{F}(t+1; \mathbb{R}^m)$
\begin{eqnarray}\label{open-loop-equilibrium-t=1}
J(t+1,X^{t+1}_{t+1}; u^{t,x,*}|_{\mathbb{T}_{t+1}})\leq J(t+1, X^{t+1}_{t+1}; (u_{t+1}, u^{t,x,*}|_{\mathbb{T}_{t+2}})).
\end{eqnarray}
Continuing the above procedure of obtaining
(\ref{system-1-t=0})-(\ref{open-loop-equilibrium-t=1}), we have for
any $k\in \mathbb{T}$
\begin{eqnarray*}\label{system-1-t=t}
\left\{\begin{array}{l}
X^{k}_{k+1} =\big{[}(A_{k,k}+\bar{A}_{k,k})X^{k}_k+(B_{k,k}+\bar{B}_{k,k})u^{t,x,*}_k+f_{k,k}\big{]}\\[1mm]
\hphantom{X^{k}_{k+1} =}+\big{[}(C_{k,k}+\bar{C}_{k,k})X^{k}_k+(D_{k,k}+\bar{D}_{k,k})u^{t,x,*}_k+d_{k,k}\big{)}w_k, \\[1mm]
X^{k}_{k} = X_k^{k-1},
\end{array}\right.
\end{eqnarray*}
and for any $u_k\in L^2_\mathcal{F}(k; \mathbb{R}^m)$
\begin{eqnarray*}\label{open-loop-equilibrium-t=t}
J(k,X^{k}_{k}; u^{t,x,*}|_{\mathbb{T}_{k}})\leq J(k, X^{k}_{k}; (u_k, u^{t,x,*}|_{\mathbb{T}_{k+1}})).
\end{eqnarray*}
Denote $\{x, X_{t+1}^{t+1}, X_{t+2}^{t+2},\cdots, X_{N-1}^{N-1},
X_{N}^{N}\}$ by  $\{x, X_{t+1}^{t,x,*}, X_{t+2}^{t,x,*},\cdots,
X_{N-1}^{t,x,*}, X_N^{t,x,*}\}\triangleq X^{t,x,*}$. Then,
$(X^{t,x,*}, u^{t,x,*})$ is an open-loop equilibrium pair.
This proves the theorem. \endpf

\subsection*{C. Proof of Lemma \ref{Lemma-Z}}\label{appendix-C}         

\emph{Proof}. Let $u^{t,x,*}_\ell=\Psi_\ell
X^{t,x,*}_\ell +\alpha_\ell, \ell\in \mathbb{T}_{k}$. Then, we have
\begin{eqnarray*}
&&\hspace{-1.5em}X^{k,t,x}_{N}=A_{k,N-1}X^{k,t,x}_{N-1}+\bar{A}_{k,N-1}\mathbb{E}_{k}X^{k,t,x}_{N-1}
+B_{k,N-1}\Psi_{N-1}X^{t,x,*}_{N-1}+\bar{B}_{k,N-1}\Psi_k\mathbb{E}_{k}X^{t,x,*}_{N-1}\\[1mm]
&&\hspace{-1.5em}\hphantom{X^{k,t,x}_{N}=}+\mathcal{B}_{N-2,N-1}\alpha_{N-1}+f_{N-2,N-1}+\big{\{}C_{k,N-1}X^{k,t,x}_{N-1}+\bar{C}_{k,N-1}\mathbb{E}_{k}X^{k,t,x}_{N-1}
\\[1mm]
&&\hspace{-1.5em}\hphantom{X^{k,t,x}_{N}=}+D_{k,N-1}\Psi_{N-1}X^{t,x,*}_{N-1}+\bar{D}_{k,N-1}\Psi_{N-1}\mathbb{E}_{k}X^{t,x,*}_{N-1}+\mathcal{D}_{k,N-1}\alpha_{N-1}+d_{k,N-1}\big{\}}w_{N-1}.
\end{eqnarray*}
To calculate $Z^{k,t,x}_{N-1}$, we need some preparations. Noting that
\begin{eqnarray*}
Z_{N}^{k,t,x}=G_{k}X^{k,t,x}_{N}+\bar{G}_{k}\mathbb{E}_{k}X^{k,t,x}_{N}+g_{k},
\end{eqnarray*}
we get
\begin{eqnarray*}
&&\hspace{-1.5em}A^T_{k,N-1}\mathbb{E}_{N-1}Z_N^{k,t,x}= A_{k,N-1}^T\mathbb{E}_{N-1}\big{[}G_{k}X_N^{k, t,x}+\bar{G}_{k}\mathbb{E}_{k}X_N^{k,t,x}+g_{k} \big{]}\\[1mm]
&&\hspace{-1.5em}\hphantom{A^T_{k,N-1}\mathbb{E}_{N-1}Z_N^{k,t,x}}=A^T_{k,N-1}G_{k}A_{k,N-1}X^{k,t,x}_{N-1}+\big{[}A^T_{k,N-1}G_{k}\bar{A}_{k,N-1}+A^T_{k,N-1}\bar{G}_{k}\mathcal{A}_{k,N-1}\big{]}\mathbb{E}_{k}X^{k,t,x}_{N-1}\\[1mm]
&&\hspace{-1.5em}\hphantom{A^T_{k,N-1}\mathbb{E}_{N-1}Z_N^{k,t,x}=}- A^T_{k,N-1}G_{k} B_{k,N-1}\Psi_{N-1}X^{t,x,*}_{N-1}\\[1mm]
&&\hspace{-1.5em}\hphantom{A^T_{k,N-1}\mathbb{E}_{N-1}Z_N^{k,t,x}=}-A^T_{k,N-1}\big{[}G_{k} \bar{B}_{k,N-1}+\bar{G}_{k}\mathcal{B}_{k,N-1}\big{]}\Psi_{N-1}\mathbb{E}_{k}X^{t,x,*}_{N-1}\\[1mm]
&&\hspace{-1.5em}\hphantom{A^T_{k,N-1}\mathbb{E}_{N-1}Z_N^{k,t,x}=}-A^T_{k,N-1}\mathcal{G}_{k}\big{[}\mathcal{B}_{k,N-1}\alpha_{N-1}-f_{k,N-1}\big{]}+A^T_{k,N-1}g_{k}.
\end{eqnarray*}
Similarly, we have the expressions of
$\bar{A}^T_{k,N-1}\mathbb{E}_{k}Z_N^{k,t,x}$,
$C^T_{k,N-1}\mathbb{E}_{N-1}\big{(}Z_N^{k,t,x}w_{N-1}\big{)}$ and
$\bar{C}^T_{k,N-1}\mathbb{E}_{k}\big{(}Z_N^{k,t,x}w_{N-1}\big{)}$.
Furthermore,
\begin{eqnarray*}
&&\hspace{-2em}Z_{N-1}^{k,t,x}=\big{[}Q_{k,N-1}+A^T_{k,N-1}G_{k}A_{k,N-1}+C^T_{k,N-1}G_{k}C_{k,N-1}\big{]}X^{k,t,x}_{N-1}\\[1mm]
&&\hspace{-2em}\hphantom{Z_{N-1}^{k,t,x}=}+\big{[}\bar{Q}_{k,N-1}+A^T_{k,N-1}G_{k}\bar{A}_{k,N-1}+A^T_{k,N-1}\bar{G}_{k}\mathcal{A}_{k,N-1}\\[1mm]
&&\hspace{-2em}\hphantom{Z_{N-1}^{k,t,x}=}+C^T_{k,N-1}G_{k}\bar{C}_{k,N-1}+\bar{A}^T_{k,N-1}\mathcal{G}_{k}\mathcal{A}_{k,N-1}+\bar{C}^T_{k,N-1}G_{k}\mathcal{C}_{k,N-1}\big{]}\mathbb{E}_{k}X^{k,t,x}_{N-1}\\[1mm]
&&\hspace{-2em}\hphantom{Z_{N-1}^{k,t,x}=}+\big{[}A^T_{k,N-1}G_{k} B_{k,N-1}+C^T_{k,N-1}G_{k} D_{k,N-1}\big{]}\Psi_{N-1}X^{t,x,*}_{N-1}\\[1mm]
&&\hspace{-2em}\hphantom{Z_{N-1}^{k,t,x}=}+\big{\{}A^T_{k,N-1}\big{[}G_{k} \bar{B}_{k,N-1}+\bar{G}_{k}\mathcal{B}_{k,N-1}\big{]}+\bar{A}^T_{k,N-1}\mathcal{G}_{k}\mathcal{B}_{k,N-1}\\[1mm]
&&\hspace{-2em}\hphantom{Z_{N-1}^{k,t,x}=}+C^T_{k,N-1}G_{k} \bar{D}_{k,N-1}+\bar{C}^T_{k,N-1}G_{k} \mathcal{D}_{k,N-1} \big{\}}\Psi_{N-1}\mathbb{E}_{k}X^{t,x,*}_{N-1}\\[1mm]
&&\hspace{-2em}\hphantom{Z_{N-1}^{k,t,x}=}+\mathcal{A}^T_{k,N-1}\mathcal{G}_{k}\big{[}\mathcal{B}_{k,N-1} \alpha_{N-1}+f_{k,N-1}\big{]}+\mathcal{{C}}^T_{k,N-1}{G}_{k}\big{[}\mathcal{D}_{k,N-1} \alpha_{N-1}+d_{k,N-1}\big{]}\\[1mm]
&&\hspace{-2em}\hphantom{Z_{N-1}^{k,t,x}=}+\mathcal{A}^T_{k,N-1}g_{k}+q_{k,N-1}\\[1mm]
&&\hspace{-2em}\hphantom{Z_{N-1}^{k,t,x}}=P_{k,N-1}X^{k,t,x}_{N-1}+\bar{P}_{k,N-1}\mathbb{E}_kX^{k,t,x}_{N-1}+T_{k,N-1}X^{t,x,*}_{N-1}+\bar{T}_{k,N-1}\mathbb{E}_kX^{t,x,*}_{N-1}+\pi_{k,N-1}.
\end{eqnarray*}
We now calculate $Z^{k,t,x}_{N-2}$. Note that
\begin{eqnarray*}
&&\hspace{-2em}A^T_{k,N-2}\mathbb{E}_{N-2}Z^{k,t,x}_{N-1}=A^T_{k,N-2}\big{[}P_{k,N-1}\mathbb{E}_{N-2}X^{k,t,x}_{N-1}+\bar{P}_{k,{N-1}}\mathbb{E}_kX^{k,t,x}_{N-1}\\[1mm]
&&\hspace{-2em}\hphantom{A^T_{k,N-2}\mathbb{E}_{N-2}Z^{k,t,x}_{N-1}=}+T_{k,{N-1}}\mathbb{E}_{N-2}X^{t,x,*}_{N-1}+\bar{T}_{k,{N-1}}\mathbb{E}_kX^{t,x,*}_{N-1}+\pi_{k,{N-1}}\big{]}\\[1mm]
&&\hspace{-2em}\hphantom{A^T_{k,N-2}\mathbb{E}_{N-2}Z^{k,t,x}_{N-1}}=A^T_{k,N-2}P_{k,N-1}A_{k,N-2}X_{N-2}^{k,t,x}\\[1mm]
&&\hspace{-2em}\hphantom{A^T_{k,N-2}\mathbb{E}_{N-2}Z^{k,t,x}_{N-1}=}+\big{(}A^T_{k,N-2}P_{k,N-1}\bar{A}_{k,N-2}+A^T_{k,N-2}\bar{P}_{k,N-1}\mathcal{A}_{k,N-2}\big{)}\mathbb{E}_kX_{N-2}^{k,t,x}\\[1mm]
&&\hspace{-2em}\hphantom{A^T_{k,N-2}\mathbb{E}_{N-2}Z^{k,t,x}_{N-1}=}+\big{[}A^T_{k,N-2}P_{k,N-1}B_{k,N-2}\Psi_{N-2}\\[1mm]
&&\hspace{-2em}\hphantom{A^T_{k,N-2}\mathbb{E}_{N-2}Z^{k,t,x}_{N-1}=}+A^T_{k,N-2}T_{k,N-1}\big{(}\mathcal{A}_{N-2,N-2}+\mathcal{B}_{N-2,N-2}\Psi_{N-2} \big{)}\big{]}X^{t,x,*}_{N-2}\\[1mm]
&&\hspace{-2em}\hphantom{A^T_{k,N-2}\mathbb{E}_{N-2}Z^{k,t,x}_{N-1}=}+\big{[}A^T_{k,N-2}P_{k,N-1}\bar{B}_{k,N-2}\Psi_{N-2}+A^T_{k,N-2}\bar{P}_{k,N-1}\mathcal{{B}}_{k,N-2}\Psi_{N-2}\\[1mm]
&&\hspace{-2em}\hphantom{A^T_{k,N-2}\mathbb{E}_{N-2}Z^{k,t,x}_{N-1}=}+A^T_{k,N-2}\bar{T}_{k,N-1}\big{(}\mathcal{A}_{N-2,N-2}+\mathcal{B}_{N-2,N-2}\Psi_{N-2} \big{)}\big{]}\\[1mm]
&&\hspace{-2em}\hphantom{A^T_{k,N-2}\mathbb{E}_{N-2}Z^{k,t,x}_{N-1}=}\times\mathbb{E}_kX^{t,x,*}_{N-2}+A^T_{k,N-2}\mathcal{P}_{k,N-1}\big{(}\mathcal{B}_{k,N-2}\alpha_{N-2}+f_{k,N-2}\big{)}\\[1mm]
&&\hspace{-2em}\hphantom{A^T_{k,N-2}\mathbb{E}_{N-2}Z^{k,t,x}_{N-1}=}+A^T_{k,N-2}\mathcal{T}_{k,N-1}\big{(}\mathcal{B}_{N-2,N-2}
\alpha_{N-2}+f_{N-2,N-2}\big{)}+A^T_{k,N-2}\pi_{k,N-1},
\end{eqnarray*}
and similar expressions for
$C^T_{k,N-2}\mathbb{E}_{N-2}\big{(}Z^{k,t,x}_{N-1}w_{N-2}\big{)}$,
$\bar{A}^T_{k,N-2}\mathbb{E}_kZ^{k,t,x}_{N-1}$
and
$\bar{C}^T_{k,N-2}\mathbb{E}_{k}\big{(}Z^{k,t,x}_{N-1}w_{N-2}\big{)}$.
Then, from (\ref{system-adjoint}) we have
\begin{eqnarray*}
&&\hspace{-2em}Z^{k,t,x}_{N-2}=\big{(}Q_{k,N-2}+A^T_{k,N-2}P_{k,N-1}A_{k,N-2}+C^T_{k,N-2}P_{k,N-1}C_{k,N-2}\big{)}X_{N-2}^{k,t,x}\\[1mm]
&&\hspace{-2em}\hphantom{Z^{k,t,x}_{N-2}=}+\big{(}\bar{Q}_{k,N-2}+A^T_{k,N-2}P_{k,N-1}\bar{A}_{k,N-2}+A^T_{k,N-2}\bar{P}_{k,N-1}\mathcal{A}_{k,N-2}+\bar{A}^T_{k,N-2}\mathcal{P}_{k,N-1}\mathcal{A}_{k,N-2}
\\[1mm]
&&\hspace{-2em}\hphantom{Z^{k,t,x}_{N-2}=}+C^T_{k,N-2}P_{k,N-1}\bar{C}_{k,N-2}+\bar{C}^T_{k,N-2}P_{k,N-1}\mathcal{C}_{k,N-2}\big{)}\times\mathbb{E}_kX_{N-2}^{k,t,x}+\big{[}A^T_{k,N-2}P_{k,N-1}B_{k,N-2}\Psi_{N-2}\\[1mm]
&&\hspace{-2em}\hphantom{Z^{k,t,x}_{N-2}=}+A^T_{k,N-2}T_{k,N-1}\big{(}\mathcal{A}_{N-2,N-2}+\mathcal{B}_{N-2,N-2}\Psi_{N-2} \big{)}+C^T_{k,N-2}P_{k,N-1}D_{k,N-2}\Psi_{N-1}\\[1mm]
&&\hspace{-2em}\hphantom{Z^{k,t,x}_{N-2}=}+C^T_{k,N-2}T_{k,N-1} \big{(}\mathcal{C}_{N-2,N-2}+\mathcal{D}_{N-2,N-2}\Psi_{N-2} \big{)}\big{]}X^{t,x,*}_{N-2}+\big{[}A^T_{k,N-2}P_{k,N-1}\bar{B}_{k,N-2}\Psi_{N-2}\\[1mm]
&&\hspace{-2em}\hphantom{Z^{k,t,x}_{N-2}=}+A^T_{k,N-2}\bar{P}_{k,N-1}\mathcal{{B}}_{k,N-2}\Psi_{N-2}+A^T_{k,N-2}\bar{T}_{k,N-1}\big{(}\mathcal{A}_{N-2,N-2}+\mathcal{B}_{N-2,N-2}\Psi_{N-2} \big{)}\\[1mm]
&&\hspace{-2em}\hphantom{Z^{k,t,x}_{N-2}=}+C^T_{k,N-2}P_{k,N-1}\bar{D}_{k,N-2}\Psi^{N-2}+\bar{A}^T_{k,N-2}\mathcal{P}_{k,N-1}\mathcal{B}_{k,N-2}\Psi_{N-2}+\bar{A}^T_{k,N-2}\mathcal{T}_{k,N-1}\\[1mm]
&&\hspace{-2em}\hphantom{Z^{k,t,x}_{N-2}=}\times\big{(}\mathcal{A}_{N-2,N-2}+\mathcal{B}_{N-2,N-2}\Psi_{N-2} \big{)}+\bar{C}^T_{k,N-2}P_{k,N-1}\mathcal{D}_{k,N-2}\Psi_{N-1}+\bar{C}^T_{k,N-2}T_{k,N-1}\\[1mm]
&&\hspace{-2em}\hphantom{Z^{k,t,x}_{N-2}=}\times\big{(}\mathcal{C}_{N-2,N-2}+\mathcal{D}_{N-2,N-2}\Psi_{N-2} \big{)}\big{]}\mathbb{E}_kX^{t,x,*}_{N-2}+A^T_{k,N-2}\mathcal{P}_{k,N-1}\big{(}\mathcal{B}_{k,N-2}\alpha_{N-2}+f_{k,N-2}\big{)}\\[1mm]
&&\hspace{-2em}\hphantom{Z^{k,t,x}_{N-2}=}+A^T_{k,N-2}\mathcal{T}_{k,N-1}\big{(}\mathcal{B}_{N-2,N-2} \alpha_{N-2}+f_{N-2,N-2}\big{)}+\bar{A}^T_{k,N-2}\mathcal{P}_{k,N-1}\big{(}\mathcal{B}_{k,N-2}\alpha_{N-2}+f_{k,N-2}\big{)}\\[1mm]
&&\hspace{-2em}\hphantom{Z^{k,t,x}_{N-2}=}+\bar{A}^T_{k,N-2}\mathcal{T}_{k,N-1}\big{(}\mathcal{B}_{N-2,N-2} \alpha_{N-2}+f_{N-2,N-2}\big{)}+C^T_{k,N-2}P_{k,N-1}\big{(}\mathcal{D}_{k,N-2}\alpha_{N-2}+d_{k,N-2}\big{)}\\[1mm]
&&\hspace{-2em}\hphantom{Z^{k,t,x}_{N-2}=}+C^T_{k,N-2}T_{k,N-1}\big{(}\mathcal{D}_{N-2,N-2}\alpha_{N-2} +d_{N-2,N-2}\big{)}+\bar{C}^T_{k,N-2}P_{k,N-1}\big{(}\mathcal{D}_{k,N-2}\alpha_{N-2}+d_{k,N-2}\big{)}\\[1mm]
&&\hspace{-2em}\hphantom{Z^{k,t,x}_{N-2}=}+\bar{C}^T_{k,N-2}T_{k,N-1}\big{(}\mathcal{D}_{N-2,N-2}\alpha_{N-2} +d_{N-2,N-2}\big{)}+\mathcal{A}^T_{k,N-2}\pi_{k,N-2}+q_{k,N-2}\\[1mm]
&&\hspace{-2em}\hphantom{Z^{k,t,x}_{N-2}}=
P_{k,N-2}X^{k,t,x}_{N-2}+\bar{P}_{k,N-2}\mathbb{E}_kX^{k,t,x}_{N-2}+T_{k,N-2}X^{t,x,*}_{N-2}+\bar{T}_{k,N-2}\mathbb{E}_kX^{t,x,*}_{N-2}+\pi_{k,N-2}.
\end{eqnarray*}
By deduction, we achieve the conclusion. This completes the proof.\endpf




\end{document}